\documentclass[pdflatex,sn-mathphys-num]{sn-jnl}

\usepackage{graphicx}%
\usepackage{multirow}%
\usepackage{amsmath,amssymb,amsfonts}%
\usepackage{amsthm}%
\usepackage{mathrsfs}%
\usepackage[title]{appendix}%
\usepackage{xcolor}%
\usepackage{textcomp}%
\usepackage{manyfoot}%
\usepackage{booktabs}%
\usepackage{algorithm}%
\usepackage{algorithmicx}%
\usepackage{algpseudocode}%
\usepackage{listings}%

\newcommand{\etype}[1]{\renewcommand{\labelenumi}{(#1{enumi})}}
\newtheorem*{prop*}{Proposition}

\def\Ann{\operatorname{Ann}}

\def\mcE{\mathcal E}

\def\tmcE{\tilde{\mathcal E}}
\def\mcI{\mathcal I}
\def\mcA{\mathcal A}

\def\mcT{\mathcal T}

\def\Af{A^{K}}

\def\pp{N,M,{\operatorname{prim}}}

\newtheorem{prop}{Proposition}[section]

\newtheorem{thm}[prop]{Theorem}
\newtheorem{conj}[prop]{Conjecture}
\newtheorem{ques}[prop]{Question}
\newtheorem{cor}[prop]{Corollary}
\newtheorem{lemma}[prop]{Lemma}

\theoremstyle{definition}
\newtheorem{Def}[prop]{Definition}
\newtheorem*{Def*}{Definition}
\newtheorem{Defs}[prop]{Definition}

\newtheorem{Facts}[prop]{Facts}
\newtheorem{Defsnot}[prop]{Definitions and notation}

\newtheorem*{notation*}{Notation}
\newtheorem*{remark*}{Remark}

\newcommand{\rad}{\operatorname{rad}}

\newcommand{\ga}{\alpha}
\newcommand{\gb}{\beta}
\newcommand{\gc}{\gamma}

\newcommand{\gve}{\varepsilon}
\newcommand{\gl}{\lambda}

\newcommand{\gvp}{\varphi}

\newcommand{\gt}{\tau}

\newcommand{\charc}{{\rm char}}

\newcommand{\lan}{\langle}
\newcommand{\ran}{\rangle}

\def\eroman{\etype{\roman}}

\newcommand{\tr}{{\rm tr}}
\newcommand{\tor}{{\rm tor}}

\newcommand{\half}{\frac{1}{2}}

\numberwithin{equation}{section}

\theoremstyle{thmstyleone}
\newtheorem{example}[prop]{Example}%
\newtheorem{remark}[prop]{Remark}%

\theoremstyle{thmstylethree}%

\begin{document}

\title[Axial identities]{Axial identities}


\author {\fnm{Louis Halle} \sur{Rowen}}
\email{rowen@math.biu.ac.il}

\affil {\orgdiv{Mathematics Department}, \orgname{Bar-Ilan University}, \\
\orgaddress{\street{Max and Anna Webb~St.}, \city{Ramat Gan}, \postcode{5290002}, \country{Israel}}}


\abstract{Axial algebras are algebras generated by semisimple idempotents whose eigenspace decompositions satisfy prescribed fusion rules, such as Jordan fusion rules, which hold in Jordan algebras and Matsuo algebras. Axial algebras are studied here in terms
 of identities involving idempotents and/or axes, and at times involving Frobenius forms. After laying out the general theory, organized as a PI-type theory of marked varieties, we turn to Jordan type, providing insight into a major theorem of I.~Gorshkov, S.~Shpectorov, and A.~Staroletov about solid subalgebras.

 This approach also leads to generic constructions of primitive $N$-generated axial algebras of Jordan type $\gl$, for $N\le 4$, and for arbitrary $N$ we obtain corresponding ``$M$-primitive'' constructions. We discuss the relevance to the problem of finding a primitive nonsingular axial algebra of Jordan type~$\half$, which is neither Jordan nor a homomorphic image of a Matsuo algebra, although we cannot provide an answer at this stage.}


%


 \pacs[Orcid ]{0000-0003-4540-9574}
\maketitle


\section{Introduction}\label{sec1}

Axial algebras --- algebras generated by semisimple idempotents --- were studied long ago by Peirce (in the associative case) and by Albert \cite{A2,A3,A4}. Over the past two decades, axial algebras have become a more active area of research since the discovery of their relation with vertex operator algebras, 3-transposition groups and the classification of finite simple nonabelian groups, nonlinear PDEs, and vector flows on manifolds, as described in the excellent overview by McInroy~\cite{M}.

It has long been known (\cite{A1}, with different terminology) that Jordan algebras generated by idempotents are axial algebras of Jordan type~$\half$, although Matsuo~\cite{Ma} found another very important class.
Recently attention has focused on describing the structure of primitive axial algebras of Jordan type~$\half$ in terms of 2-generated axial subalgebras (which are easily seen to be Jordan algebras).
Gorshkov, Shpectorov, and Staroletov \cite{GSS} and Desmet~\cite{D} proved results describing when 2-generated subalgebras of a primitive axial algebra of Jordan type $\half$ over a field $F$ are ``solid,'' meaning all primitive idempotents $\ne 0,1$ are axes.

This paper was motivated by the desire to understand these impressive results in terms of a PI-like theory involving ``idempotental'' (resp.~``axial'') identities of axial algebras, in other words, identities in which some of the substitutions are reserved for idempotents (resp.~axes). In a sense, this is a special case of ``generalized identities,'' which define a variety, thereby permitting tools of PI-theory paralleling \cite[Chapter~7]{Row}, and generic objects. Idempotental identities can be found for Jordan algebras, under the name \textbf{Peirce derivatives} in \cite{LCS}, and as \textbf{Peirce $s$-identities} in \cite{McC}; McCrimmon uses them to reprove quickly that the Albert algebra is exceptional. However, these are more subtle than ordinary polynomial identities, as seen in the discussion of ``marked'' varieties.

We indicate how the main theorems of \cite{D} and \cite{GSS} could be interpreted through axial identities (although the basic ideas remain the same), determining when idempotents are axes satisfying the Jordan fusion rules. At the outset, we outline a heuristic computer-free argument in characteristic not dividing 30, for the theorem of \cite{GSt} that 3-generated primitive axial algebras of Jordan type $\half$ are Jordan algebras.

In a different direction,
we construct generic (also called ``universal'') axial algebras of Jordan type $\gl$, before directly imposing a Frobenius form.\footnote{\cite{HRS1} already constructed a version of universal axial algebras with a Frobenius form, but the approach here is directly from the viewpoint of identities and algebraic varieties.} In order to develop this properly, we need to consider axial algebras over a commutative Noetherian ring $C$ (rather than just a field).

Conceptually, we start within the theory of universal algebra; identities can be formulated as elementary universal sentences, where idempotents are viewed as distinguished elements.
An interesting complication arises from the condition that an idempotent be ``primitive,'' which states that $ax = x$ implies $x\in Ca;$ this cannot be formulated as a universal sentence. Thus when imposing an $\forall\exists$ sentence on a marked variety, we do not have a single variety with a generic object, but rather obtain a finite number of marked varieties (for finite-dimensional algebras). We circumvent this difficulty by introducing an indeterminate for the coefficient of $a$ in $ax.$

This connects with the issue of finite-dimensionality of $N$-generated primitive axial algebras of Jordan type $\half$, which we circumvent by considering words of bounded length.\footnote{In fact, using the classification of finite simple groups, Cuypers and Hall \cite{CuH} proved that finitely generated 3-transposition groups are finite, which solves the problem for Matsuo algebras.} For $N\ge 4$, we obtain $N$-generated generic $M$-primitive axial algebras of Jordan type~$\half$ over integral domains, for suitably large~$M$. For $N=4$ \cite{DRS} proved finite-dimensionality, so we obtain all the nonisomorphic 4-generated generic primitive axial algebras of Jordan type~$\half$ over integral domains. Some of them are Matsuo, some are Jordan, and the remaining ones, \emph{if they exist}, are neither Jordan nor a homomorphic image of a Matsuo algebra. The same sort of analysis works for nonsingular axial algebras.

 In the process, we deduce from the other axioms that the $0$-eigenspace of a primitive axis of Jordan type $\gl$ is a subalgebra.

As with \cite{D}, our focus is on the algebraic structure via idempotents. Although we concentrate on the commutative case, the noncommutative case also is of interest and is mentioned briefly in Definition~\ref{not2} and in~\S\ref{dg}.

To summarize, we present
\begin{itemize}
    \item Structural identities: formulation of axial/idempotental/Frobenius identities and their consequences.

\item Solid-subalgebra application: reinterpretation of the Gorshkov–Shpectorov–Staroletov/Desmet machinery.

\item
Generic constructions: finite families of generic (M)-primitive objects, with genuine primitive objects when axial admissibility is known.
\end{itemize}

\section{Fundamental notions}
\begin{Defsnot}\label{not1}

In this paper, $C$ always denotes a commutative Noetherian ring containing the element $\half$. An element $c\in C$ is called \textbf{regular} if, for each $c'\in C,$ $c c' =0$ implies $c' =0.$ 
In the literature, $C$ would be
 an infinite field~$F$, but the more general situation is useful for handling generic algebras.

 We review some familiar definitions.
\begin{enumerate}

\item
For $a\in A,$ write $L_a$ for the left multiplication map
$z\mapsto az$, and $R_a$ for the right multiplication map
$z\mapsto za$.

 \item Even when $C$ is not a field, we say that $A_\gl$ is a left $\gl$-\textbf{eigenspace} of~$a$ if $$A_\gl = \{y\in A: ay = \gl y\}.$$
 \item
 A \textbf{left axis}
$a$ is an 
idempotent for which $A$ is a direct sum of its
left eigenspaces with respect to $L_a$. A \textbf{right axis} is
defined analogously. We write $L_a^k(y)$ for $a(a(\dots(ay))),$ taken $k$ times.

 \item A \textbf{(2-sided) axis} $a$ is a left axis
which is also a right axis, for which $L_a R_a = R_a L_a$. Thus we can write $aya$ without ambiguity.

 In the commutative case, which we assume throughout except in~\S\ref{dg}, $L_a= R_a$, so every left axis is an axis. The eigenvalues $1$ (and possibly $0$) play a special role, so from now on we designate $S(a)$ as the eigenvalues of $a$ other than $ 0,1$.

 \emph{
We assume throughout that $\gl$, $1-\gl$, and $\gl -\gl'$ are invertible in $C$ for each distinct $\gl, \gl'\in S(a).$ In particular, when $S(a) = \{\gl\},$ we assume that $\gl (1-\gl)$ is invertible in $C$.} (One could have $S(a)=\emptyset$.)
 Then, for any axis~$a,$ $A$~decomposes as the direct sum of $C$-subspaces
\begin{equation}
 \label{dec17} A= A_{0,1}(a) \oplus \left(\bigoplus_{\gl \in S} A_{\gl}(a)\right),
\end{equation}
where $A_\gl(a)$ denotes the $\gl$-eigenspace of~$L_a$, and $A_{0,1}(a) := A_1(a) \oplus A_0(a) $.
Often we write $S $ for $ S(a)$, and call $a$ an $S$-\textbf{axis}; so the eigenvalues of $a$ are in $S\cup \{0,1\}.$ Thus, any element $y\in A$ is written uniquely as
\begin{equation}
 \label{dec2}
y = y_1+ y_0 +\sum _{\gl \in S} y_{\gl} , \qquad ay_1 =y_1, \quad ay_0 = 0, \quad ay_\gl= \gl y_\gl. 
\end{equation}

\item For $y\in A$, define the \textbf{annihilator} $\Ann_C y = \{c\in C : cy =0\}.$

\item Obviously $a\in A_1(a) $ for any idempotent $a$. The idempotent $a$ is \textbf{primitive} if $A_1(a) = C a;$ when $C$ is an integral domain with $\Ann_C(a)=0$, we write $y_1 = \gvp_a(y) a$ for $\gvp_a(y)\in C.$

 \item An algebra $A$ is \textbf{axial over} $\mathcal E$ if $A$ is generated by a set $\mathcal E$ of axes.\footnote{The generating set of axes often is denoted in the literature as $X$, but we want to reserve the letter $X$ for indeterminates.} Usually, we suppress $\mathcal E$ in the notation. An axial algebra $A$ is $N$-\textbf{generated} if $|\mcE|=N;$
 in this case we write $\mcE = \{ a_i : 1\le i \le N\}.$
 $A$~is $S$-\textbf{axial} if each axis in $\mathcal E$ is an $S$-axis.
 We focus on the case $|S|=1,$ but \textbf{Monster type} is included in the case $|S|=2$.

\item The \textbf{basic fusion rules} are the two rules:
\begin{enumerate}
 \item $A_{0,1}(a)$ is a $C$-subalgebra of $A$;
\item
 $A_{\gl} (a)$ is
 an $A_{0,1}(a)$-module for each $\gl\in S$.
\end{enumerate}

\item An $S$-axis $a$ is of \textbf{Jordan type $S$} if $a$ satisfies the basic fusion rules and the rules $$A_{\gl}(a)A_{\gl'}(a) \subseteq \delta_{\gl,\gl'} \,A_{0,1}(a), \qquad \forall \gl,\gl'\in~S,$$ where $\delta$ denotes the Kronecker delta.

\item An \textbf{axial algebra of Jordan type $S$}, called an \textbf{axial algebra of Jordan type~$\gl$} when $S = \{\gl\},$ is an $S$-axial algebra for which each $a\in \mathcal{E} $ is of Jordan type $S$.\footnote{In the literature, Jordan type includes the fusion rule $A_0(a)^2 \subseteq A_0(a)$, but we shall see that this comes automatically for primitive axes of Jordan type $\gl$ over a field.}

\item A \textbf{PAJ-$S$} is a primitive axial algebra of Jordan type $S$. 
A \textbf{PAJ-$\gl$} is a PAJ-$\{\gl\}.$

\item A \textbf{marked homomorphism} of axial algebras is a homomorphism of nonassociative algebras which sends the generating axes to axes.

\noindent The remaining items are used because we do not know yet whether all $N$-generated PAJ-$\half$ over a field are finite-dimensional.
\item For an axial algebra $A$ over $\mathcal E$, define $A^{(M)}:=A^{(M)}(\mcE)$ to be the $C$-submodule spanned by all words in $\mcE$ of length $\le M$ (permitting any placement of parentheses). A class of axial algebras $\mcA$ is \textbf{strongly finite} if there is some $M$, which we call $\Phi(N),$ for which
$A = A^{(M)}$ for every $A\in \mcA.$ The number $N$ is $\gl$-\textbf{axially admissible} if the class of $N$-generated PAJ-$\gl$ is strongly finite.

\item An axis $a$ is $M$-\textbf{primitive} if $ay=y$ for $y\in A^{(M)}$ implies $ y\in Ca.$

\item An \textbf{MPAJ-$\gl$} is a (commutative) $\gl$-axial algebra over an integral domain for which each axis $a\in \mathcal{E} $ is $M$-primitive of Jordan type $\gl$.
\end{enumerate}
\end{Defsnot}

\begin{remark}\label{sp1}
 If $|\mcE|=N$ then $A^{(M)}$ is spanned by some number $\theta(N,M)\le (4N)^M$ elements, seen by counting all products of elements of $\mcE$.
\end{remark}

\subsubsection{Facts from commutative ring theory}

For later reference, we recall some basic facts which impact on axial algebras.

\begin{remark}\label{quo} A commutative ring $C$ is \textbf{reduced} if $c^2 = 0$ for $c\in C$ implies $c = 0.$
\begin{enumerate}\eroman
\item $C$ has a maximal nilpotent ideal $\mathcal N:= \operatorname{Nil}(C),$ and $\bar C : =C/\mathcal N$ is reduced.

 \item $\bar C $ has a finite set of minimal prime ideals $P_1,\dots, P_k$ whose intersection is 0.

 \item
 The quotient ring $Q(\bar C )$ obtained by inverting regular elements of~$\bar C $ is isomorphic to the direct product $Q(\bar C/P_1)\times \dots \times Q(\bar C/P_k)$ of fields.

 \item For $K=Q(\bar C )$, and $A$ a torsion-free $\bar C$-algebra,
 define $\Af = A\otimes_{\bar C } K.$
 \begin{enumerate}
 \item Any element in $\Af$ can be written $c^{-1}y$ for regular $c\in \bar C ,$ $y\in A.$

 \item Any idempotent of $\Af$ has the form $c^{-1}a $ where $a^2 = ca.$ Conversely, if $a ^2 = c a$ regular for $c\in \bar C$, then $ c^{-1}a$ is idempotent.
 \end{enumerate}
 \end{enumerate}
\end{remark}

\subsubsection{Miyamoto maps}

There is a way of obtaining new axes.

\begin{Defs}\label{Miy7}
 \begin{enumerate}
 \item Suppose that $a$ is an $S$-axis of $A$. Picking nonempty $S' \subseteq S$, define the \textbf{Miyamoto map} $\gt_{a,S'} :A \to A$ by
 $y \mapsto y_{0,1}+ \sum _{\gl \in S\setminus S'} {y_{\gl}}\, - \sum _{\gl \in S'} y_\gl$. $\gt_{a,S'}$ clearly is a $C$-linear map of $A$, whose square is the identity map.

\item A \textbf{Miyamoto involution} is a Miyamoto map which is an
 algebra automorphism.
 \item A \textbf{Miyamoto axis} is an axis $a$ whose Miyamoto maps all are Miyamoto involutions. 
\end{enumerate}
\end{Defs}

As in \cite[Definition~4.1]{HSS2}, we note:

\begin{lemma}
 Any axis of Jordan type $S$, not necessarily primitive, is a Miyamoto axis.
\end{lemma}
\begin{proof} An easy observation on 2-graded algebras, noting that
 for any $S' \subseteq S$, $A$ is 2-graded by~\eqref{dec17}, where $A_{\gl}(a)$ is odd for $\gl \in S'$ and all other eigenspaces of $a$ are even.
\end{proof}

 \subsubsection{Frobenius forms}
 \begin{Defs}
 \begin{enumerate}
 \item An \textbf{associative bilinear form} $(\phantom{w},\phantom{w}): A \times A \to C$ satisfies
 $$\left(\sum \ga_{i }x_i,\sum \gb_{ j}y_j\right) = \sum _{i,j} \ga_i \gb_j (x_i,y_j), \quad (xy,z) = (x,yz), \qquad \forall x,y,z\in A.$$

 \item A \textbf{Frobenius form} is a symmetric associative bilinear form.\footnote{Any associative bilinear form on a commutative axial algebra is symmetric, cf.~\cite[Proposition~3.5]{HRS1}.}

\item An axis $a$ is \textbf{anisotropic} if $(a,a)$ is invertible in $C$.\footnote{Note the invertibility condition, which restricts to nonzero over a field.} A Frobenius form is \textbf{anisotropic} if every primitive axis is anisotropic. A Frobenius form is \textbf{normal} if $(a,a)=1$ for every primitive axis $a.$

\item The \textbf{radical} $\rad A$ of a Frobenius form on $A$ is $\{x\in A: (x,A)=0\}.$ 
 \end{enumerate}
 \end{Defs}
\begin{lemma} \label{normal} $ $\begin{enumerate}
 \item \cite[Theorem~2.3]{RoSe6} Every PAJ-$\gl$ $A$ over a field has a normal Frobenius form, given by $(a,y) = \gvp_a(y)$ for every primitive axis $a$ and all $y\in A$.
 \footnote{This also holds more generally for many PAJ-$S$, by \cite{RoSe6}.}

 \item Suppose $a$ is an axis of an $S$-axial algebra $A$ with an anisotropic Frobenius form, with $y,y' \in A$. Notation as in \eqref{dec2},
 \begin{enumerate}
 \item $(y_\gl,y'_{\gl'})=0$ for all $\gl\ne \gl'\in S$.

 \item $(y_\gl,y'_{0})= (y_\gl,a)= 0$ for all $\gl\in S$.

 \item $(a,y_0)=0.$

 \item If $a$ is primitive then $(a,a)y_1 = (a,y)a$. Thus, for $a$ anisotropic, $y_1 = \frac{(a,y)}{(a,a)}a$.

 \item If $A$ is spanned by a set $\tilde{\mathcal E}$ of primitive axes, then
 the following are equivalent:
 \begin{enumerate}
 \item $y \in \rad A$.

 \item $(a,y) =0$ for every $a\in \tilde{\mathcal E}.$

 \end{enumerate}
 \end{enumerate}
\end{enumerate}
\end{lemma}
\begin{proof}

 1. It should be noted that originally the result was proved by~\cite[Theorem~4.1]{HSS1} when $A_0(a)^2 \subseteq A_0(a)$ for all primitive axes $a$.

2. (a) $\gl (y_\gl,y'_{\gl'}) = (ay_\gl,y'_{\gl'}) = (y_\gl,ay'_{\gl'}) =\gl'(y_\gl,y'_{\gl'}) ,$ and $\gl-\gl'$ is presumed invertible.

 \indent \indent (b) As in (a), since $\gl$ and $1-\gl$ are presumed invertible.

 \indent \indent (c)
 $(a,y_0) = (a^2,y_0)= (a,ay_0) =(a,0)=0.$

 \indent \indent (d) Write $y_1 = \gc a.$ Then $ \gc (a,a)a = ( a,\gc a)a = (a,y_1)a = (a,y)a,$ by (a), (b), and (c).

 \indent \indent (e) $(ii)\Rightarrow (i)$ Write $z = \sum \gc_i a_i. $ Then $(z,y) = \sum \gc_i (a_i,y) = 0.$

 \indent \indent \indent \ $(i) \Rightarrow (ii)$ Obvious.
\end{proof}

\begin{remark}\label{rady}
\begin{enumerate}\eroman
\item \cite{HRS1} classifies the axial algebras of ``Monster type'' that have a Frobenius
form. As noted in \cite[Conjecture 8.15]{M}, Franchi, Mainardis, and Shpectorov conjecture that every axial algebra of Monster type has a Frobenius
form.
 \item The Frobenius form passes to subalgebras, and also to direct products
 where distinct components are taken to be orthogonal,
 but need not pass to a homomorphic image since it might not be well-defined.

\item Lemma~\ref{normala} permits an alternate description of the radical of a $\gl$-axial algebra, without using the Frobenius form. When $A$ is spanned by a set $\tmcE$ of primitive axes, we say $A$ is \textbf{nonsingular} if for each $0\ne y\in A$ there is $a\in \tmcE$ for which $a(ay) \ne \gl ay.$
\end{enumerate}
 \end{remark}
\subsubsection{Annihilators and torsion}

Suppose $A$ is an algebra over $C.$
Coefficients from $C$ can interfere with the theory, by annihilating elements of $A.$ To circumvent this difficulty, we introduce another notion.

\begin{Def}\label{annid}
 \begin{enumerate}
 \item For a subset $\mcI \subseteq A$, $\Ann_C (\mcI) = \cap _{f\in \mcI} \Ann_C(f).$

 \item $\Ann_A(c) = \{a\in A: ca= 0\}.$

 \item $\mcT$ will always denote a lattice of ideals of $A$, i.e., closed under sums and intersections, such that $\Ann_A c \in \mcT$ and $c\mcI \in \mcT$ for every $c\in C$ and $\mcI \in \mcT.$

 \item A $\mcT$-\textbf{annihilator} in $C$ is an annihilator of some nonzero $\mcI \in \mcT$.
 \end{enumerate}
\end{Def}

We recall some standard facts about $\mcT$-annihilators.

\begin{lemma}
 \label{tor1} Suppose $A$ is an algebra over $C.$
 \begin{enumerate}\eroman
 \item If $P$ is a $\mcT$-annihilator in $C$ then $P =\Ann_C (\Ann_A P)$.

 \item Every maximal $\mcT$-annihilator in $C$ is a prime ideal.

 \end{enumerate}
\end{lemma}
\begin{proof} (i) $P \subseteq \Ann_C (\Ann_A P)$ is obvious. Likewise if $P = \Ann _C(\mcI),$ then $\mcI \subseteq \Ann _A(P),$ so $P \supseteq \Ann_C (\Ann_A P)$.

(ii) Suppose $P$ is an $\mcT$-annihilator of $\mcI\in \mcT$. If $c_1,c_2 \in C$ such that $c_1c_2\in P,$ then $c_1c_2 \mcI =0.$ If $c_2\notin P$ and $c_2\mcI \ne 0 $, then $\Ann_C (c_2 \mcI)$ contains both $P$ and $c_1,$ contrary to the maximality of $P$ unless $c_1\in P$.
\end{proof}

 \subsection{ The case $S = \{\gl\}$}

 We explore the situation of \cite{HRS2}, namely for $S = \{\gl\}$, where $\gl \ne 0,1$.
 To simplify the exposition, we assume that $(1-\gl)^{-1} \gl^{-1}\in C.$

\begin{example}\label{prim0} Let $a$ be a primitive axis of Jordan type $\gl$, and write $y = y_1 + y_0 + y_\gl $ as in \eqref{dec2}.

 \begin{equation}\label{ay}
 ay = y_1 + \gl y_\gl.
\end{equation}
\begin{equation}\label{ay2}
 a(ay) = y_1 + \gl^2 y_\gl,
\end{equation}
\begin{equation}\label{ay3}
 a(a(ay)) = y_1 + \gl^3 y_\gl,
\end{equation}
\begin{equation}\label{ay21}
 a(ay)-\gl ay = (1-\gl)y_1,
\end{equation}
 which lets us solve for $y_1$, and
\begin{equation}\label{ay321}
 a(ay)- ay = (\gl-1)\gl y_\gl,
\end{equation}
 which lets us solve for $y_\gl$.

These can be encoded in the equations

\begin{equation}\label{ax10}
 a(a(ay)-\gl ay) = a(ay-\gl y) ,
\end{equation}
\begin{equation}
 \label{ax10a} a(a(ay)- ay )= \gl( a(ay)- ay),
\end{equation}
 both of which are equivalent to

\begin{equation}
 \label{ax1} a(a(ay)) = (\gl+1)( a(ay))- \gl ay.
\end{equation}
 \end{example}

 \begin{lemma}\label{axfree}
 Conversely, suppose an idempotent $a\in A$ satisfies \eqref{ax1}. Then
 \begin{enumerate}\eroman
 \item $a$ is an $\gl$-axis.

 \item $A_1(a)$ is spanned as a $C$-module by all expressions of the form $$a(ay -\gl y), \qquad y\in A.$$

 \item $A_\gl(a)$ is spanned as a $C$-module by all expressions of the form $$\frac{1}{\gl}(a(ay) -ay), \qquad y\in A.$$

\item $A_0(a)$ is spanned as a $C$-module by all expressions of the form \begin{equation}
\label{a0} y+\frac{1}\gl a(ay) -\frac{\gl+1}\gl ay.
\end{equation}
 \end{enumerate}
 \end{lemma}

\begin{proof} Recall that \eqref{ax1} implies both
\eqref{ax10} and \eqref{ax10a},
which
respectively imply $ a(ay-\gl y) \in A_1(a)$ and $a(ay)- ay\in A_\gl(a)$. Then
noting that
\begin{equation}\label{ax10b}
 (\gl-1)y_1 = \gl ay -a(ay), \qquad (\gl-1)y_\gl = \frac{1}{\gl}(a(ay)-ay)
\end{equation}
we see that \begin{equation}
 (\gl-1)( y_1 + y_\gl) = - (\gl-1)\frac{1}\gl a(ay) +\frac{\gl^2-1}\gl ay;
\end{equation}

canceling $(\gl-1)$
shows that $y_0= y-(y_1+y_\gl)$ is of the form
 \eqref{a0}.
\end{proof}

This explicit description does not depend on the base ring $C.$

\emph{Assume for the rest of this subsection that $A$ has a normal Frobenius form}, and $a$ is a primitive axis\footnote{The results hold mutatis mutandis for an anisotropic Frobenius form.}. Then
 \begin{equation}\label{ayf}
 ay = (a,y)a + \gl y_\gl.
\end{equation}
\begin{equation}\label{ayf2}
 a(ay) = (a,y)a + \gl^2 y_\gl,
\end{equation}
\begin{equation}\label{ayf21}
 a(ay) = \gl ay + (1-\gl)(a,y)a.
\end{equation}
To have $y_\gl \in A_\gl(a)$ we need
\begin{equation}
 \label{prim01} ay - (a,y)a =\gl y_\gl = a y_\gl = \frac 1 \gl a(ay - (a,y)a) = \frac 1 \gl a(ay) - \frac 1 \gl(a,y)a .
\end{equation}

We shall also need, from \eqref{ayf21},
\begin{equation}\label{ayf31}\begin{aligned}
 a( a&(ay))= \gl a(ay) + (1-\gl)(a,y)a \\& = \gl ^2 ay + \gl(1-\gl)(a,y)a +(1-\gl)(a,y)a = \gl ^2 ay + (1-\gl^2)(a,y)a ,
\end{aligned}
\end{equation}
which encodes the fact that the minimal polynomial of $L_a$ divides $x(x-1)(x-\gl).$

Furthermore, \eqref{prim01} becomes $ ay - (1-\frac{1}{\gl})(a,y)a = \frac{1}{\gl}aya,$ or
\begin{equation}\label{aya100} \gl ay + (1-{\gl})(a,y)a =a(ay).\end{equation}

\begin{lemma}\label{typ}
 Every element of $A_{0}(a)$ is of the form $
\gl y +(1-\gl )(a,y)a -ay.$
Every element of $A_{0,1}(a)$ is of the form $
\gl y-ay.$
\end{lemma}
\begin{proof} For the first assertion, one checks directly that $a(\gl y +(1-\gl )(a,y)a -ay) = 0.$
 Conversely if $y\in
A_0(a)$ then $ (a,y)=0$, so $$y= \gl^{-1}(\gl y + (1-\gl )(a,y)a -ay).$$

The second assertion follows easily.
\end{proof}

\begin{lemma}
 \label{normala}
 Suppose $A$ is a PAJ-$\gl$ with an anisotropic Frobenius form, and $y \in A$.
 Then $y\in \rad A$ if and only if $a(ay) = \gl ay $ for every $a\in \tilde{\mathcal E}.$
\end{lemma}
\begin{proof}
 Recall from Lemma~\ref{normal}(2)(e) that $y\in \rad A$ if and only if $(a,y) =0,$ i.e., $y=y_0 +y_\gl.$

 $(\Rightarrow)$ If $(a,y)=0$ then $ay =ay_0 +ay_\gl = \gl y_\gl,$ so $a(ay) = \gl^2 y_\gl =\gl ay.$

 $(\Leftarrow)$ $y_1=0$ by Lemma~\ref{axfree}(ii).
\end{proof}

\begin{lemma}
 \label{aya}
 Suppose $a,b$ are primitive axes of Jordan type $\gl$. Then
 \begin{enumerate}\eroman
 \item $aba -bab = (1-\gl)(a,b)(a-b)$.
 \item $ab = 0$ if and only if $(a,b)=0$ and $b_\gl = 0.$
 \end{enumerate}
\end{lemma}
\begin{proof}
(i) By \eqref{ayf21} and its reverse for $b,$
 \begin{equation}
 aba = a(ab) =\gl ab + (1-\gl)(a,b)a,
 \end{equation}
 \begin{equation}
 bab = \gl ba + (1-\gl)(a,b)b,
 \end{equation}
 so subtract.

 (ii) By \eqref{ayf}.
\end{proof}





\subsubsection{Miyamoto maps for a PAJ-$\gl$}

For $S=\{\gl\}$, we write
 $\gt_a :=\gt_{a,\gl} : y \mapsto y^{\gt_a } = y -2 y_\gl.$

\begin{lemma}
 \label{Miy5} (As in \cite[Proposition~3.6]{DRS})
\begin{enumerate}\label{Miy4}\eroman
 \item If $a,b$ are primitive Miyamoto axes then $b^{\gt_a}$ is a Miyamoto axis, and $2ab$ is a linear combination of Miyamoto axes.
\item If $A$ has a Frobenius form with a primitive $\gl$-axis $a $ in which $(a,a)=1$, and $y\in A$, then
 \begin{enumerate}
 \item $ay = \frac{\gl}{2} (y - y^{\gt_a}) +(a,y)a$.

 \item
$ y^{\gt_a} = y + \frac{2}{\gl}(a, y)a- \frac{2}{\gl}ay.$
 \end{enumerate}
 \end{enumerate}
\end{lemma}
\begin{proof}

 (i) The first assertion holds since $\gt_a$ is an algebra automorphism, and the second assertion follows since $ay = {\gl} y_\gl +\varphi_a(y) a = \frac{\gl}{2} (y - y^{\gt_a}) +\varphi_a(y) a$.

 (ii) $ay = {\gl} y_\gl +(a,y)a = \frac{\gl}{2} (y - y^{\gt_a}) +(a,y)a$.
 We solve to get (b).
\end{proof}

\begin{prop}\label{Miy6} Suppose $A$ is an axial algebra over a set of Miyamoto axes $\mathcal E$.
 \begin{enumerate} \eroman
 \item $A$ is spanned by a set $\tmcE$ of Miyamoto axes obtained from $\mathcal E$.

 \item If $y = \sum \gc_i x_i $, then $y_\mu = \sum \gc_i (x_i)_\mu$ for any eigenvalue $\mu$ of $a.$
 \end{enumerate}
\end{prop}
\begin{proof}
 (i) Apply induction to Lemma~\ref{Miy5}(i).

 (ii) Decompose into eigenspaces of $a.$
\end{proof}

Proposition~\ref{Miy6} is critical in many proofs. Even for axial algebras of Monster type, one can still use 2-gradings to obtain Miyamoto involutions, as noted in \cite[Proposition~3.4]{HRS1} and \cite[Theorem~2.15]{W}.

\subsubsection{2-generated PAJ-$\gl$}\label{3gen2}

For primitive $\gl$-axes $a,b$
 \cite[\S 3]{GSS} calls $A= \lan\lan a,b \ran\ran$ \textbf{toric} for $(a,b)\notin \{0,1\};$ \textbf{flat} for $(a,b)=0$; and \textbf{baric} for $(a,b)=1$.
Toric algebras have an important
parametrization for $\gl = \half:$

\begin{lemma}(\cite[Lemma~ 3.5(b)]{GSS},\cite[Lemma~3.2]{D})\label{tor}
 Any toric algebra has a basis $\{e, u, f\},$ with $u$ being the identity element, \[
\textstyle{e^2=f^2=0,\quad ef=\frac{1}{8}u,}
\] and
\begin{equation}
 (e,e)=(f,f)=(u,e)=(u,f)=0, \qquad (e,f) =\frac{1}{4},\qquad (u,u)=2.
\end{equation}

 An element $y$ is a nontrivial idempotent if and only if \begin{equation}
 \label{ids}
 y = \gve e + \gve^{-1} f + \half u, \quad \gve \in F^\times.
 \end{equation}
\end{lemma}

 \subsubsection{Matsuo algebras}

 Matsuo algebras comprise the most important class of axial algebras after Jordan algebras.
\begin{Defs} \cite[Example~2.7]{M}
 A \textbf{Matsuo algebra of type~$\gl$} is an axial algebra of Jordan type~$\gl,$ such that for any two distinct axes $a,b$ either $(\gt_a\gt_b)^2$ is the identity map, in which case $ab=0,$ or \begin{equation}
 \label{Mats0} ab = \frac{\gl}{2}\left(a+b-b^{\gt_a}\right).
 \end{equation}
\end{Defs}

\begin{lemma}\label{M1} \cite{Ma}, \cite{HRS2,HSS1}
 Every axis of a Matsuo algebra $A$ of type $\gl$ is primitive, and $A$ is endowed with a normal Frobenius form satisfying $(a,b)= \frac{\gl}{2}$ for any non-orthogonal axes $a\ne b.$
\end{lemma}
(We recall that in the case of \eqref{Mats0},
 $b_1 + \gl b_\gl = ab = \frac{\gl}{2}(a+b-b^{\gt_a}) = \frac{\gl}{2}a+ \gl b_\gl,$ implying $b_1 = \frac{\gl}{2}a.$ Thus $a$ is primitive, and $A$ has a normal Frobenius form, where by inspection $(a,b) = \frac{\gl}{2}.$)

\begin{lemma}\label{Matl}
 Over a field (as in \cite[Lemma~2.11]{D}), \begin{enumerate}\eroman
 \item $ab=0$ or $(a,b) = \frac{\gl}{2},$ for any axes $a\ne b$ of a Matsuo algebra of type~$\gl.$
 \item Conversely, if $ab=0$ or $(a,b) = \frac{\gl}{2},$ for any axes $a\ne b$ of a PAJ, then $A$ is Matsuo of type $\gl$.

 \item $A$ is Matsuo of type~$\gl$ if and only if
 \begin{equation}\label{Mat2}
\left(a(ab)-\gl(ab)-\frac{\gl(1-\gl)}{2}a\right)(ab) = 0
 \end{equation}
 for any axes $a\ne b$ of $A$.

 \item $A$ is Matsuo of type~$\gl$ if and only if \begin{equation}\label{Mat2a}
\left(aba-bab-\frac{\gl(1-\gl)}{2}(a-b)\right)(ab) = 0
 \end{equation} and
 \begin{equation}\label{Mat3a}
\left(a(ab)-\gl(ab)-\frac{\gl(1-\gl)}{2}a\right)(ab-b) = 0
 \end{equation}
are both satisfied, for all axes $a,b$ of $A.$
 \end{enumerate}
\end{lemma}
\begin{proof}
 (i) By Lemma~\ref{M1}.

(ii) Reverse the argument of the first line of the proof of Lemma~\ref{M1}.

(iii) Recall that $ a(ab)-\gl ab = (1-\gl)(a,b)a$, by \eqref{ayf21}. Thus
\begin{equation}
 \label{ay4} a(ab)-\gl(ab)-\frac{\gl(1-\gl)}{2}a = (1-\gl)\left((a,b)-\frac{\gl}{2} \right)a.
\end{equation}

$(\Rightarrow)$ In view of \eqref{ay4}, if $(a,b)= \frac{\gl}{2} $ then $ a(ab)-\gl ab- (1-\gl)(a,b)a =0.$

$(\Leftarrow)$
$0= (1-\gl)\left((a,b)-\frac{\gl}{2} \right)(ab),$ by \eqref{ay4}.

Thus either $(a,b)-\frac{\gl}{2} $ is 0, implying $(a,b) = \frac{\gl}{2} ,$ or otherwise
$0 = a(ab) = (a,b)a+\gl^2 b_\gl,$
 implying $(a,b)=0$ and $b_\gl = 0$, so $ab =0$ by Lemma~\ref{aya}(ii). Hence we are done by (ii).

 (iv) $(\Rightarrow)$ As in (iii).

 $(\Leftarrow)$ Assume that $(a,b)\ne\frac{\gl}{2}$. \eqref{Mat2a} says $ (1-\gl)\left((a,b)-\frac{\gl}{2}\right)(a-b)(ab) = 0 $, so $0 = (a-b)(ab) = aba-bab = (1-\gl)(a,b)(a-b),$
 by Lemma~\ref{aya}(i), implying $a= b$ or $(a,b) = 0.$
 Likewise \eqref{Mat3a} says $(a,b)=\frac{\gl}{2}$ or $a(ab) = ab,$ which implies $b_\gl = 0.$ Thus, if $a\ne b$ then $ab=0$ by Lemma~\ref{aya}(ii). \end{proof}

\begin{remark} The following facts were already known. \label{MatJ} \begin{enumerate}\eroman
 \item
By \cite[Theorem 7.1]{HSS2} together with \cite[Theorem~6.4]{HRS2}, any PAJ-$\gl$ for $\gl\ne \half$ is a homomorphic image of a Matsuo algebra of type~$\gl$.

 \item Examples of Matsuo algebras of Jordan type $\half$ that are not Jordan algebras are given in \cite{DMR},\cite{DPSC}, \cite{W}, and \cite{GSS}.

 \item In \cite{DMR,GSS}, Jordan algebras are given which are not homomorphic images of Matsuo algebras.
\end{enumerate}
\end{remark}

The \textbf{projective axial graph} of a PAJ spanned by $\tmcE$ is the graph whose vertices are the axes in
$\tmcE,$ with an edge between $a$ and $b$
 if $\gvp_a(b)\ne 0,$ i.e., if
$(a,b) \ne 0$ when there is an anisotropic Frobenius form.

An axial algebra $A$ is \textbf{connected} if its projective axial graph is connected.
\footnote{This is stronger than saying $ab \ne 0$, as evidenced by the 3-dimensional flat algebra.}

 One major classification question of axial algebras has been:

 \begin{ques}\label{Spe}
 (Shpectorov) Are connected PAJ-$\half$ either Jordan or homomorphic images of
 Matsuo algebras?
 \end{ques}

 \subsection{Identities}

Intuitively, in universal algebra, an identity of an algebra is a universally quantified equation holding for all substitutions in the algebra.
An \textbf{(algebraic) variety} is a class of algebras defined by identities. Varieties are characterized by the properties that they are closed under products, subalgebras, and homomorphic images.

\subsubsection{Jordan identities and power-associativity}

The most celebrated variety of commutative nonassociative algebras is the class of Jordan algebras, defined by identities of degree~4. We also glance at the noncommutative version.

\begin{Defs}\label{not2}

\begin{enumerate}
\item A \textbf{flexible algebra} satisfies the identity $(xy)x = x(yx).$

\item A \textbf{noncommutative Jordan algebra} is a flexible algebra satisfying the identity \begin{equation}
 \label{nonc} ((xx)y)x = (xx)(yx),
\end{equation}

called
the ``Jordan identity''.

\item A \textbf{Jordan algebra} is a commutative algebra
which satisfies the Jordan identity.

\item An \textbf{almost Jordan} algebra is a commutative algebra satisfying the following consequence \cite{O2} of the Jordan identity:\begin{equation}
 \label{alm}
 2((yx)x)x + y((xx)x) = 3(y(xx))x,\quad\forall x,y\in A. \footnote{This terminology comes from the title of \cite{HP}. Jacobson~\cite[p.~307ff]{J} defined a \textbf{Lie triple product} in terms of the equation $$L_{(xy)z-x(yz)}=[L_z,[L_x,L_y]], \qquad \forall x,y,z \in A,$$ which in \cite[eq.~5, p.~1115]{O2} is seen by linearization to be equivalent to \eqref{alm}, and is called a \textbf{Lie triple algebra}, which is more general than a Lie triple system. The situation is further complicated by another equivalent definition used in \cite{O3}. Hence \cite{AL,HP,O2,O3,P1,P2,Si1,Si2} are all about the same notion.}
\end{equation}

\item Specializing $y$ to $x$ in the Jordan identity yields $((xx)x)x = (xx)(xx),$ which is called
\textbf{4-power-associativity}.

\end{enumerate}
 \end{Defs}


\subsubsection{Background results}

Here are some results concerning identities of an arbitrary commutative algebra $A$ over a field $F$, including the Jordan identity and 4-power-associativity.

\begin{Facts}\label{br}
 \begin{enumerate}\eroman


 \item Any commutative 2-generated axial algebra of Jordan type $\half$ is a Jordan algebra. (Verifying the Jordan identity is an easy exercise, to be discussed below.)

 \item \cite[Corollary to Theorem 4]{O1} If $A$ has a unit element and satisfies an identity of degree $\le 4$ not implied by commutativity, then $A$ either satisfies~(4) or (5) of Definition~\ref{not2}, or the identity \begin{equation}\label{th}
 2(x(x y^2)) + 2(x^2y)y +(xy)^2 = 2((xy)y)x + 2((xy)x)y+x^2y^2.
 \end{equation}
(\cite[Corollary to Theorem 5]{O1} classifies all commutative algebras with unit satisfying an identity of degree $\le 5$.)

 \item \cite{CaHP} classifies all commutative algebras satisfying an identity of degree $\le 4$.

 \item \cite[p.~1114, l.~-12, without proof]{O2} When $\frac 13 \in F,$ a commutative almost Jordan algebra is Jordan if and only if it satisfies 4-power-associativity.

 \item \cite[Equation (6)]{O2} Any idempotent of an almost Jordan algebra is a $\half$-axis.

 \item \cite[Lemma~1]{O2} Any idempotent of an almost Jordan algebra is of Jordan type $\half$. This is a key fact, since, in conjunction with~(v), it implies that any axial algebra which is almost Jordan is of Jordan type $\half.$


 \item \cite[Theorem 1]{O2} An almost Jordan algebra $A$ is Jordan if and only if $A_1(a)$ and $A_0(a)$ are Jordan algebras for some axis $a$. In particular, if $A$ is simple, then it is a Jordan algebra. (This is improved in (viii).)

\item \cite[Theorem 1]{Si1} (stated for Lie triple algebras, which was noted to be the same as almost Jordan algebras) Let $\mathcal J(A)$ denote the ideal of an almost Jordan algebra $A$ generated by $((xx)x)x - (xx)(xx),$ for all $x\in A.$ Then $\mathcal J(A)^2 = 0,$ and $A/\mathcal J(A)$ is Jordan. In particular, an almost Jordan 4-power-associative algebra is Jordan.

 \item \cite[ p. 1115, Equations (4) and (5)]{O2} $A$ is almost Jordan if and
only if $D_{x,y}: = [L_x,L_y]$ is a derivation for all $x, y \in A.$

 \item \cite[Lemma~4, p.~554]{A2} If the characteristic of $F$ is prime to 30, then 4-power-associativity of a commutative algebra implies power-associativity (i.e., all powers of $x$ associate, for any $x\in A$).\footnote{\cite{Ko} showed that 4-power-associativity and the identity $X^3X^2 = X^4X$ imply power-associativity in characteristic 3; 4-power-associativity and the identity $X^4X^2 = X^5X$ imply power-associativity in characteristic 5.}

\item \cite[Theorem 2, p.~559]{A2} in characteristic prime to 30, and \cite[Theorem 3]{Ko} in characteristics 3 and~5: In any commutative power-associative algebra, all idempotents are $\half$-axes satisfying the basic fusion rules, as well as some weaker fusion rules.\footnote{$A_\half A_i \subseteq A_\half + A_{1-i}$ and $A_i^2 \subseteq A_i$, for $i=0,1$.}

 \item \cite[Proposition~6.4]{D} Any almost Jordan algebra spanned by idempotents is Jordan.

 \item \cite[Proposition~A8]{ChG} Suppose the characteristic of $F$ is prime to~30, and $A$ satisfies an identity of degree $\le 4$ satisfied by some algebra with a unit element, but not implied by
commutativity. If $A$ has a nondegenerate Frobenius form then $A$ is a Jordan algebra.\footnote{Their method is elegant, going over the list of \cite{O1}. To remove the requirement of a unit element, one could take instead the list of \cite{CaHP}, which is restated in the list of \cite[Theorem~1]{RB}.
(1) of their list reduces to 4-power-associativity for any algebra with an idempotent. (3) is disposed of in \cite{AL}. \cite{RB} proves that (5) implies the identity $xy^2 = (xy)y,$ and taking $y=a$ idempotent, Example~\ref{prim0} yields $\gl^3 = \gl^2,$ which we exclude. This leaves (2), which is $$ 2\gb ((xy)^2 - x^2y^2) + \gc (((xy)x)y + ((xy)y)x - (y^2x)x - ( x^2y)y ) = 0.$$}

 \item The class of $N$-generated PAJ-$\gl$ is strongly finite, for $N\le 4$. More precisely,
 $A= A^{(2)}$ for $N=2,$ by \cite{HRS2}, $A= A^{(3)}$ for $N=3,$ by \cite{GSt}, and $A= A^{(7)}$ for $N=4,$ by \cite{DRS}\footnote{For $N=4$ it suffices to take words of the form $L_{a_1}\dots L_{a_6}(a_7)$.}. So far the situation for $N>4$ is unclear.

 \end{enumerate}
\end{Facts}

\begin{remark}
\begin{enumerate}\eroman
 \item By item (v), a PAJ-$\gl$ for $\gl\ne \half$ cannot be almost Jordan, and its theory differs in many ways from a PAJ-$\half$.

 \item Items (viii), (ix), (x), and (xiii) do not require any idempotents!

 \item In view of (xiii), any PAJ-$\half$ with unity, whose Frobenius form is nondegenerate, satisfying \eqref{th}, is a Jordan algebra.
 \end{enumerate}
\end{remark}

\begin{ques} Which degree 5 identities of \cite[Corollary to Theorem~5]{O1} are satisfied by every PAJ-$\half$?
\end{ques}
\subsubsection{Digression: Noncommutative analogs}\label{dg}

We digress briefly to discuss the more intricate noncommutative situation. Some classes of noncommutative axial algebras of Jordan type are given in \cite{RoSe5, RoSe6}, together with the normal Frobenius form obtained in \cite{RoSe6}.

Next, some low degree identities. The exterior (Grassmann) algebra satisfies the identity $[[X_1,X_2],X_3]=0$ which is of degree 3, but is not relevant here since none of its degree 1 elements are idempotent.

 The left and right alternative identities, respectively $x^2y = x(xy)$ and $xy^2 = (xy)y$, have degree~3. (Compare with Fact~\ref{br}(ii), which is for commutative algebras.) Also $2\times 2$ matrix algebras over a field satisfy a standard identity; this is essentially the only associative case generated by two idempotents, as seen in \cite{RoSe1}.

 \cite{Si2} introduces (in Russian) a noncommutative version of a Lie triple algebra.

 \begin{Facts}
 \begin{enumerate}\eroman
 \item (Albert~\cite[Theorem 1, p. 574]{A2}) A flexible algebra is noncommutative Jordan if and only if $L_x, R_x, L_{x^2}, R_{x^2}$ all commute, for each $x.$

 \item (Albert~\cite[p. 553]{A2}) Given any algebra $A$ and $\mu\in C,$ there is an algebra $A(\mu)$ with the same $C$-module structure as $A$ but new multiplication $x\circ y =\mu xy + (1-\mu)yx.$ Any axis of $A$ remains an axis of $A(\mu)$ (since it is obviously an idempotent, and the 2-sided eigenvalue decomposition remains). Thus $A(\mu)$ is an axial algebra if $A$ is an axial algebra.

 \item Write $A^+$ for $A(\half)$ in (ii). $A^+$ is a commutative axial algebra when $A$ is axial, and is Jordan when $A$ is noncommutative Jordan.

\item Albert~\cite[Theorem 5, p. 562]{A2} recovers most of Fact~\ref{br}(xi) in the noncommutative case.\footnote{ (But the proof of Fact~\ref{br}(ix) relies heavily on commutativity.) Albert \cite[p.~563]{A2} also gives an example of a flexible power-associative algebra with a nontrivial idempotent which does not satisfy the Jordan fusion rules.}
 \end{enumerate}
 \end{Facts}

\section{Initial results}

\subsection{Jordan type implies $A_0(a)^2 \subseteq A_0(a)$}

First we resolve the discrepancy between our definition of Jordan type and the customary definition in the literature.

\begin{thm}\label{3gen}
 $A_0(a)^2 \subseteq A_0(a)$ for any anisotropic primitive axis $a$ satisfying the basic fusion rules, in an algebra $A$ with a Frobenius form.
\end{thm}
\begin{proof}
 If $y,z \in A_0(a)$ then $yz \in A_{0,1}(a)$ by the basic fusion rule. Writing $yz = d_0 +\ga a$ for $d_0 \in A_0(a),$ we have $(a,a)\ga = (a, d_0 +\ga a)= (a,yz) = (ay,z) = (0,z) = 0,$
 so $yz=d_0 \in A_0(a).$
\end{proof}

\begin{cor}\label{3gen1}
 $A_0(a)^2 \subseteq A_0(a)$ for any PAJ $A$ of Jordan type $\gl$ over a field.
\end{cor}
 \begin{proof}
$A$ has a normal Frobenius form by \cite[Theorem~2.3]{RoSe6}.
(\cite{HSS1} yields a Frobenius form, but uses the property $A_0(a)^2 \subseteq A_0(a)$, which we are trying to prove.)
 \end{proof}

Note that \cite[Theorem~2.3]{RoSe6} gives a much broader result providing a normal Frobenius form, in view of \cite[Theorem A]{RoSe5}.

\subsection{3-generated PAJ-$\gl$}\label{3gen3}

\begin{remark} Gorshkov and Staroletov
 \cite{GSt} proved the important theorem that if $A = \lan\lan a,b,c \ran \ran$ for primitive axes $a,b,c,$ of Jordan type $\gl,$ then $A$ is spanned by $\{ a,b,c, ab, ac, bc, (ab)c, (ac)b , a(bc)\}$, and they also provided an explicit multiplication table, described concisely in \cite{GSS}. 

 \cite{DRS} provides a technique for proving the Gorshkov-Staroletov theorem, which relies only on the normal Frobenius form and Miyamoto involutions, focusing on 4-generated axial algebras, but
 does not explicitly state a proof of the Gorshkov-Staroletov theorem, so we review the argument briefly.
They work with Miyamoto automorphisms, defining $[[a_1, \dots, a_\ell]] = \gt_{a_1}\dots \gt_{a_\ell}$.
 \eqref{ayf21} is their equation (1.1), used to prove \cite[Proposition 1.2 and Proposition 2.1(iv)]{DRS} which can be seen to yield the desired result.
 \end{remark}

\section{Idempotental and axial identities}

Let us modify the notion of identities, along the lines of generalized identities taken from the theory of associative algebras, to serve in the context of axial algebras.

\begin{Defs}\label{idm} All words are nonassociative.
 \begin{enumerate}

 \item
A \textbf{generalized monomial} (over $A$) is a word $$f(y_1,\dots, y_m;X_1,\dots, X_n),$$ where the $y_i$ are fixed elements of $A,$ for example $y_1 (X_1(y_2X_2))$ or~$(y_1 (X_1y_2))X_2$. The $X_i$-\textbf{degree} is the number of times that $X_i$ appears.
The $y_i$ are called the ``coefficients'' of the generalized monomial. For technical reasons, we do not require all of the $y_i$ to appear in the generalized monomial.

 \item A \textbf{generalized polynomial} is a formal linear combination $f(y_1,\dots, y_m;X_1,\dots, X_n)= \sum _i c_i h_i(y_1,\dots, y_m;X_1,\dots, X_n)$ of generalized monomials (so the coefficients come from each monomial).

 \item A \textbf{generalized Frobenius monomial (resp.~polynomial) } is a generalized monomial (resp.~polynomial), but with the language expanded to contain a Frobenius form $(\phantom{w},\phantom{w})$.

 \item
 An \textbf{idempotental monomial} is a word $$h(E_1,\dots,E_m; X_1,\dots, X_n),$$ written $h(E;X)$ or $h$ for short, where $E_i,X_j$ are indeterminates, and we always replace $(E_iE_i)$ by $E_i.$ In the commutative setting, we rearrange the $E_i$ in   alphabetical order. For example, we replace $E_2(E_3(E_1X_1)) $ by $((E_1X_1)E_3)E_2. $ A \textbf{word} in the $E_i$ is an idempotental monomial in which only the $E_i$ appear.

 \item An \textbf{idempotental polynomial} is a formal linear combination over $C$ of
idempotental monomials.

\item A \textbf{generalized $V$-idempotental polynomial} is defined analogously to an idempotental polynomial, but with coefficients in $A$ as in generalized polynomials.

 \item A (generalized) idempotental polynomial $f(E_1,\dots,E_m; X_1,\dots, X_n)$ is of \textbf{degree}~$d_i$ in $E_i$ if $E_i$ is of degree $\le d_i$ in each monomial of $f,$ $d_i$ minimal such; $f(E_1,\dots,E_m; X_1,\dots, X_n)$ is of \textbf{degree}~$d'_i$ in $X_i$ if $X_i$ is of degree $\le d_i'$ in each monomial of $f,$ $d'_i$ minimal such.
 The idempotental polynomial $f$ is \textbf{homogeneous} in the indeterminate ~$X_j$ if each of its monomials is of the same degree in ~$X_j$.

 When $f$ is homogeneous of degree 1 in $X_j$, we call $f$ \textbf{linear} in $X_j.$ $f$ is \textbf{pure} if $f$ has degree 0 in each $X_j.$
\end{enumerate}
\end{Defs}
We need to pass to linear generalized polynomials.

\begin{Defs}
 \begin{enumerate}
 \item The $
 X_j$-\textbf{linearization step} of~$f(E_1,\dots,E_m; X_1,\dots, X_n)$~is
 \begin{equation}
 \begin{aligned}
 \Delta_j (f) (E_1,\dots, &E_m; X_1,\dots, X_n, X_{n+1})\\ := &f(E_1,\dots, E_m; X_1,\dots,X_j + X_{n+1}, X_{j+1} ,\dots, X_n) \\ &-f(E_1,\dots,E_m ; X_1,\dots,X_j,\dots X_n)\\ &\quad -f(E_1,\dots,E_{m}; X_1,\dots,X_{j-1},X_{n+1},X_{j+1},\dots, X_n).
 \end{aligned}
 \end{equation}

 \item The $
 X_j$-\textbf{linearization} of $f$ is the application of the $X_j$-linearization step successively until one arrives at a polynomial which is linear in $X_j.$

 \item The \textbf{multilinearization} of $f$ is obtained by linearizing each indeterminate.
 \item
 $f_{E_m\mapsto X_{n+1}}(E_1,\dots,E_{m-1}; X_1,\dots, X_{n+1})$ denotes the polynomial obtained by substituting $X_{n+1}$ for $E_m$ throughout in $f$.
 \end{enumerate}
\end{Defs}

\begin{example}
\cite[eq. 6]{A3} 4-power-associativity linearized in $X_2$ yields
\begin{equation}\label{PA1}
 4(X_1X_2)X_1^2 = 2((X_1X_2)X_1)X_1 + (X_1^2X_2) X_1 + X_1^3X_2.
\end{equation}

\end{example}

 Sometimes one specializes back after multilinearizing.

\begin{Defs}
 \begin{enumerate}
 \item A \textbf{generalized identity} of $A$ (written GI) is a generalized polynomial which vanishes under all ``evaluations,'' i.e., specializations of the indeterminates $X_i$ to elements of $A$.


 \item
 A $V$-\textbf{idempotental identity} of an algebra $A$ is an idempotental polynomial $$f(E_1,\dots,E_m; X_1,\dots, X_n)$$
 such that $f(E;X)\mapsto 0$ under all idempotental specializations of the $E_i$ to idempotents of $V$, and of the indeterminates $X_j $ to elements of $ A.$ We call these \textbf{idempotental specializations}. Thus, the idempotental specializations of the $E_i$ give generalized monomials over~$A$.

 \item
 Similarly, a \textbf{primitive} $V$-\textbf{idempotental identity} of an algebra $A$ is an idempotental polynomial $$f(E_1,\dots,E_m; X_1,\dots, X_n)$$
 such that $f(E;X)\mapsto 0$ under all  specializations of the $E_i$ to primitive idempotents of $V$.


 \item
 A $V$-\textbf{axial identity} of $A$ is an idempotental polynomial such that $f(E;X)\mapsto 0$ under all specializations of the $E_i$ to primitive axes of Jordan type $S$ in~$V$ and of $X_j $ to elements of $ A.$ ($S$ is given, usually $\{\gl\}$.) The axial identity is \textbf{pure} if $f$ is pure.
\end{enumerate}
\end{Defs}

So we have the following hierarchy for a marked variety:

$\text{ordinary identity}
\subset
\text{idempotental identity}
\subset
\text{axial identity}.$

\begin{example} When $A$ is spanned by primitive idempotents,
specializing $X_2$ in \eqref{PA1} to an idempotent $E_1$ shows that 4-power-associativity is equivalent to the idempotental identity

\begin{equation}
 \label{powa} 4(X_1 E_1)X_1^2 = 2((X_1 E_1)X_1)X_1 + (X_1^2 E_1) X_1 + X_1^3 E_1.
\end{equation}\end{example}

\begin{remark} Suppose $A$ is an algebra.
 \begin{enumerate}
 \item Clearly GIs of $A$ pass to subalgebras (containing the coefficients) and homomorphic images (with respect to the images of the coefficients). If we have a product $\prod_{k\in I}A_k$ of algebras, and $f(y_{1,k},\dots, y_{m,k};X_1,\dots, X_n)$ is a GI of $A_k,$ then $f(y_{1},\dots, y_{m};X_1,\dots, X_n)$ is a GI of $\prod A_k,$ where $y_i$ is the vector $(y_{i,k}).$ In this sense, GIs define what one may call a \textbf{marked variety} where we mark the $y_i$ and their homomorphic images $\bar y_i$.
 \item Accordingly, we define an \textbf{identically defined GI} to be a GI whose coefficients $y_i$ are taken arbitrarily in a set $\mathcal S_\mathcal{I}$ defined by a set $\mathcal{I}$ of identities. Our main case is $\mathcal S_\mathcal{I}$, the set of idempotents, defined by the identity $X^2-X=0.$ The same argument as in (1) goes through when we mark the elements of $\mathcal S_\mathcal{I}$, and take the coefficients in their image but we must still be careful: For a surjection $ A \to \bar A,$ $\overline{\mathcal S_\mathcal{I}(A)} \subseteq \mathcal S_\mathcal{I}(\bar A),$ but the inclusion may be strict. For example, $\bar A$ may contain new idempotents.

 We say that all algebras satisfying a given set $\mathcal S$ of identically defined GIs define a
 \textbf{marked variety defined by $\mathcal S$}.
 \item An \textbf{idempotental variety} is a marked variety defined by idempotental identities. An \textbf{axial variety} is a marked variety defined by axial identities.

 \item The standard results about homogeneous polynomials in \cite[Corollary~2.3.30]{Row} say that
 for an infinite integral domain $C,$ marked varieties are closed under tensor extensions, with respect to the same set of marked elements.

 \item Primitivity of axes is \emph{not} preserved under direct products, thereby complicating the theory.
 \end{enumerate}
\end{remark}

 Now assume for the rest of this paper that $A$ is spanned by idempotents.

 \begin{remark} \begin{enumerate}
 \item The most important cases are for $V= C a$, for $V = \lan\lan a,b \ran\ran $ for primitive idempotents $a,b$ in $V$, to deal with problems of solid subalgebras, and for $V = A.$

 \item Although GIs pass to homomorphic images, Frobenius forms do not, cf.~Remark~\ref{rady}(ii).
 \end{enumerate}

 \end{remark}

\subsection{Axial Frobenius identities}

Primitive axes are not defined in terms of identities until we introduce a Frobenius form into the language, in which case the formulation of Lemma~\ref{normal}(2)(d) is relevant.
\begin{Defs} When $A$ is endowed with a Frobenius form,
 \begin{enumerate}
 \item An \textbf{idempotental Frobenius polynomial} is an idempotental polynomial $f(E;X)$ where the $E_i,X_j$ are as above, but with the language also containing a Frobenius form $(\phantom{w},\phantom{w})$. We include the appearances of $E_i$ and $X_j$ inside the Frobenius form when counting the degree. In particular $f(E;X) := f(E_1,\dots,E_m; X_1,\dots, X_n)$ is linear in $E_m$ if $E_m$ appears exactly once in each summand of $f,$ including within the Frobenius brackets.

 \item
 A \textbf{$V$-idempotental Frobenius identity} is an idempotental Frobenius polynomial which is an identity in the language of algebras equipped with a Frobenius form,
 in which the variables $E_i$ are specialized to idempotents in $V$.


 \item
 A \textbf{$V$-axial Frobenius identity} of $A$ is an idempotental Frobenius polynomial $f:=f(E;X)$ which vanishes when the variables $E_i$ are specialized to primitive axes of Jordan type $S$ in $V.$


\end{enumerate}
\end{Defs}

Since we are looking for identities, we need
a process $$f_{E_m\mapsto X_{n+1}}(E_1,\dots,E_{m-1}; X_1,\dots, X_{n+1})$$ to reduce the number of idempotents, passing from
 a $V$-idempotental (Frobenius) identity eventually to an $A$-(Frobenius) identity. The linear case is easy.

\begin{lemma} Suppose that $A$ is spanned by primitive axes, and $f(E_1,\dots,E_m; X_1,\dots, X_n)$ is an axial Frobenius identity.
 If $f$ is linear in $E_m,$ then $f_{E_m\mapsto X_{n+1}}$ also
 is an axial Frobenius identity. In particular, when $m=1,$ $f_{E_1\mapsto X_{n+1}}$ is an axial Frobenius identity.
\end{lemma}
\begin{proof}
 Any element $y_{n+1}$ of $A$ can be written as $\sum_k \gc_k a_{m,k}$ for primitive axes $a_{m,k},$ so
 $$f_{E_m\mapsto X_{n+1}}(a_1, \dots, a_{m-1}; y_1, \dots, y_{n+1})= \sum \gc_k f(a_1, \dots, a_{m-1},a_{m,k}; y_1, \dots, y_n),$$ which is $0.$
\end{proof}

\subsubsection{Idempotental Frobenius identities for a primitive ${\gl}$-axis}

\begin{example}\label{prim00}
To prove that an idempotent $a$ of $A$ is semisimple with eigenvalues contained in ~$\{ 0,\gl ,1\}$, we need that
$$0 =L_a(L_a-1)(L_a -\gl)(y),$$ which follows from the decomposition of $y,$ and yields the $X_1$-linear generalized identity
\begin{equation}
 \label{pi3} a(a(aX_1) - (\gl +1)aX_1 + \gl X_1).
\end{equation}

Equation \eqref{ayf21} yields the $X_1$-linear generalized Frobenius identity
\begin{equation}
 \label{aya10}(1-\gl)(a,X_1)a +\gl a X_1 -a(aX_1).
\end{equation}

 \eqref{aya10} shows that, in the presence of a normal Frobenius bilinear form, primitivity translates to the $X$-linear axial Frobenius identity
\begin{equation}\label{aya1}
 \gl E_1 X_1 +(1- \gl) (E_1,X_1)E_1 -E_1(E_1X_1).
\end{equation}

Also we need to know when $ y\in Ca + A_\gl(a)+A_0(a),$ which holds whenever $ay - (a,y)a\in A_\gl(a),$ which by \eqref{aya100} again becomes an idempotental Frobenius identity.
\end{example}

\begin{remark}
 \label{FI} Given a $C$-algebra $A,$ together with a form $( \phantom{w},\phantom{w}) :A\times A\to C,$ we can characterize $( \phantom{w},\phantom{w}) $ being a Frobenius form by means of the bilinear identities
 $( \sum_i c_iX_i, \sum_j c_j'X_j) = \sum_{i,j}c_ic_j'( X_i, X_j),$
 the ``associative'' identity
 $( X_1X_2,X_3) = (X_1,X_2X_3) $, and the   identity of symmetry $(X_1,X_2) = ( X_2, X_1) $.

 Then we get normality and primitivity of $E_i$ by applying
 $(E_i,E_i)=1$ and \eqref{aya1}.\end{remark}

\subsubsection{Axial Frobenius identities for a single primitive axis}


\begin{thm}
 \label{prim} Suppose $S = \{\gl_1,\dots, \gl_t\}$ in an integral domain $C$ in which $\gl_i$ and $\gl_i-1,$ are invertible for each $1\le i \le t,$ and $\gl_i-\gl_j$ is invertible whenever $i\ne j$. \begin{enumerate}\eroman
 \item Given an $S$-axis $a$, there is an $X_1$-linear generalized polynomial, which we designate $f_{\gl_i}(a;X_1)$, which takes values precisely in $A_{\gl_i}(a)$ for each~$1\le i\le t.$

 \item An idempotent $a$ is an $S$-axis if and only if $a$ satisfies specific $X_1$-linear idempotental identities (determined by $S$).
 \item When $A$ has a Frobenius form, an anisotropic $S$-axis~$a$ is primitive if and only if it satisfies a specific $X_1$-linear Frobenius identity $(E_1,X_1)E_1 +(E_1,E_1) g_S(E_1;X_1)$ (determined by~$S$).
 \end{enumerate}
 \end{thm}
\begin{proof}

(i) $L_a^j y = y_1 + \sum_{i=1}^{t} \gl_i^j y_{\gl_i}.$  Write these $t+1$ equations for $1\le j\le t+1$. Invert  the Vandermonde matrix $(\gl_i^j:\, 0\le i \le t,\ 1\le j\le t+1 )$, taking $\gl_0=1$, noting that the Vandermonde determinant of \cite[Proposition~2.3.26]{Row} is seen to be invertible since, up to sign, it is a product of $\gl_i$ $\gl_i -1$, and $\gl_i -\gl_j$.
We can solve by Cramer's rule for $y_{\gl_0}=y_1$ and all $y_{\gl_i}$ respectively by $X_1$-linear polynomials $f_{1} (E_1;X_1)$ and $f_{\gl_i} (E_1;X_1)$, $1\le i \le t,$ and then obtain $y_0$ via $X_1- f_1 (E_1;X_1)- \sum _{i=1}^t f_{\gl_i}(E_1;X_1) $.

(ii) $E_1 f_{\gl_i}(E_1;X_1) - \gl _i f_{\gl_i}(E_1;X_1) $ yields an $X_1$-linear  axial identity for each $\gl_i\in S$, and likewise $E_1 f_1(E_1;X_1) - f_1(E_1;X_1) $ and $E_1 f_0(E_1;X_1)$.

(iii) $(\Rightarrow)$ Streamlining (i),
 $L_a^j y = \frac{(a,y)}{(a,a)}a + \sum_{i=1}^{t} \gl_i^j y_{\gl_i},$ so we can use these $t+1$ equations to eliminate the $y_{\gl_i}$ and have some expression only involving the $L_a^j y$ and $(a,y)a,$ which thereby, multiplying through by $(a,a),$ yields an $X_1$-linear idempotental Frobenius identity.

$(\Leftarrow)$ We reverse the proof of $(\Rightarrow).$ If $(E_1,X_1)E_1 + (E_1,E_1)g_S(E_1;X_1)$ is a $Ca$-idempotental Frobenius identity, then $$ (a,y)a + (a,a)g_S(a;y) =0 $$ for any $y,$ but the same argument used in proving $(\Rightarrow)$ shows
 $$y_1 + g_S(a;y) =0 ,$$ and thus $y_1 = \frac
 {(a,y)}{(a,a)}a \in C a,$ proving $a$ is primitive.
\end{proof}

\begin{remark}
 Putting $\gc = (a,a),$ we could rewrite the idempotental Frobenius polynomial of the proof as $(E_1,X_1)E_1 + \gc g_S(E_1;X_1),$ which vanishes precisely on the marked specialization to those primitive axes $a'$ for which $(a',a')=\gc.$ This relates to \cite[Problem~3.2]{GSh}.
\end{remark}

\subsubsection{Fusion rules as axial Frobenius identities}\label{axid}

Jordan type $S = \{\gl\}$ translates to the following conditions:
 \begin{enumerate} \eroman
 \item $ A_{\gl}(a) A_{0}(a)\subseteq A_{\gl}(a).$
 \item $A_{0}(a)^2 \subseteq A_{0,1}(a).$
 \item $ A_{\gl}(a)^2 \subseteq A_{0,1}(a).$
 \end{enumerate}

\begin{remark}\label{SerLem} For a primitive axis $a,$
 \textbf{Seress' Lemma} says that (i) implies
\begin{equation}
 \label{Se} a(yz) = (ay)z + a(y_0z_0) ,
\end{equation}
for any $ y \in A$ and $z\in
A_{0,1}(a).$

 Seress' Lemma is proved over a field by noting since $a$ is primitive that $y_1 z_0 = y_0 z_1 = 0,$ and \begin{equation}
 \begin{aligned}
 a((y_1 + y_\gl + y_0)(z_1 +z_0)) & = a(y_1z_1 + y_\gl z_1 +y_0z_1 + y_1z_0 +y_\gl z_0 + y_0z_0)\\ & = y_1z_1 + \gl y_\gl z_1 +\gl y_\gl z_0 + a(y_0z_0)
 = (ay)z +a(y_0 z_0).
 \end{aligned}
 \end{equation}

 $A_{0}(a)^2\subseteq A_{0}(a)$ by Corollary~\ref{3gen1}. Hence the last term in \eqref{Se} is $0$. For more general base rings, one could appeal to Theorem~\ref{3gen} for any anisotropic primitive axis $a$ satisfying the basic fusion rules, in an algebra with a Frobenius form.
\end{remark}
\begin{lemma}\label{fusf2}
 Seress' Lemma implies Condition (i).
\end{lemma}
\begin{proof}
\eqref{Se} implies Condition (i) since $a (y_\gl z_0) = (a y_\gl ) z_0 = \gl (y_\gl z_0)$.
 \end{proof}


 This leaves us to deal with the other two conditions (ii) and (iii), which can easily be translated into axial identities (albeit more complicated) using Lemma~\ref{axfree}. We treat (iii); (ii) is handled analogously.

\begin{lemma}\label{fusf3} Condition (iii) is equivalent to $$ a((ay)(az)) - \gl ^2 (a,z)ay -\gl ^2 (a,y)az = \gl ^2(ay,z)a +(1-3\gl^2) (a,y)(a,z)a,
 $$ i.e., the axial Frobenius identity

 \begin{equation}
 \begin{aligned}
 E_1((E_1X_1) & (E_1X_2))- \gl ^2 (E_1,X_2)E_1X_1 -\gl ^2 (E_1,X_1)E_1X_2 =\\ & \gl ^2(E_1X_1,X_2)E_1 +(1-3\gl^2) (E_1,X_1)(E_1,X_2)E_1.
 \end{aligned}
 \end{equation}
\end{lemma}

\begin{proof}
$x=(ay-(a,y)a)(az-(a,z)a)=\gl^2y_\gl z_\gl.$

 Condition (iii) is equivalent to
$ax=(a,x)a.$
 Noting that $(a,ay) = (a^2,y)= (a,y)$, and also
$$(a,(a,z)a(ay)) = (a,(a,z)ay) = (a,z)(a,y)=(a,y)(a,z)=(a,(a,y)a(az) ),$$ we have

 \begin{equation}
 \begin{aligned}
 a((ay)(az))& - (a,z)a(a(ay)) -(a,y)a(a(az))\\ &=(a,(ay)(az) )a- 2(a,y)(a,z)a
 \end{aligned}
 \end{equation}

 Using \eqref{ayf21}, which also yields $$(ay,az) = (a(ay),z) = (\gl ay + (1-\gl)(a,y)a,z) =\gl(ay,z) +(1-\gl)(a,y)(a,z),$$ we obtain
 \begin{equation}\begin{aligned}
 (a,(ay)&(az)) = (a(ay),az) = (\gl ay +(1-\gl)(a,y)a,az) \\ &= \gl (\gl(ay,z) +(1-\gl)(a,y)(a,z))+(1-\gl) (a,y)(a,z).
\end{aligned}
\end{equation}

Now \eqref{ayf31} yields

\begin{equation}
 \begin{aligned}
 a((ay)(az)) & - (a,z)(\gl ^2ay + (1-\gl^2)(a,y)a)\\ &\qquad -(a,y)(\gl ^2 az + (1-\gl^2)(a,z)a)\\ &=(\gl (\gl(ay,z) +(1-\gl)(a,y)(a,z))\\ &\qquad +(1-\gl) (a,y)(a,z) )a- 2(a,y)(a,z)a ,
 \end{aligned}
 \end{equation}
 or
\begin{equation} \label{pi2}
 a((ay)(az)) - \gl ^2 (a,z)ay -\gl ^2 (a,y)az = \gl ^2(ay,z)a +(1-3\gl^2) (a,y)(a,z)a. \end{equation}

\end{proof}

\begin{remark}
 \label{summ} Since the generalized Frobenius identities given here are multilinear in the $X_i$, they can be verified by specializing the~$X_i$ to primitive axes, in view of~Proposition~\ref{Miy6}.
In particular, it is enough to check the Jordan fusion rules for an idempotent $a$ in an algebra generated by $a$ and two primitive $\gl$-axes of Jordan type.
\end{remark}

\subsubsection{Idempotental Miyamoto identities}

Likewise, there is an axial Frobenius identity to determine when Miyamoto maps are necessarily Miyamoto involutions. Lemma~\ref{Miy5}(ii) lies at the heart of both \cite{DRS} and \cite{D}.

\begin{lemma} All Miyamoto maps of a PAJ-$\gl$ $A$ with a normal Frobenius form are
 Miyamoto involutions if and only if $A$ satisfies the axial Frobenius identity \begin{equation}
 \begin{aligned}
 &(E_1, X_1X_2)E_1- E_1(X_1X_2) \\ & = (E_1, X_1)E_1X_2 +(E_1, X_2)E_1X_1 - (E_1X_1)X_2 -(E_1X_2)X_1 \\ & \qquad +\frac{2}{\gl}((E_1, X_1)(E_1, X_2)E_1+ (E_1X_1)(E_1X_2)) \\ & \qquad \qquad - \frac{2}{\gl} ((E_1, X_1)E_1(E_1X_2) + (E_1, X_2)E_1(E_1X_1)).
 \end{aligned}
 \end{equation}
\end{lemma}
\begin{proof} In view of Lemma~\ref{Miy5}(ii)(b),
\begin{equation}
 \label{Mi1} \begin{aligned}
 y^{\gt_a} z^{\gt_a}& = \left(y -\frac{2}{\gl}ay + \frac{2}{\gl}(a,y)a\right)\left(z -\frac{2}{\gl}az + \frac{2}{\gl}(a,z)a\right) \\& = yz - \frac{2}{\gl} y(az) + \frac{2}{\gl}(a,z)ay \\& -\frac{2}{\gl} (ay)z +\frac{4}{\gl^2}(ay)(az) -\frac{4}{\gl^2} (a,z)a(ay)
 \\& +\frac{2}{\gl}(a,y)az -\frac{4}{\gl^2}(a,y)a(az)+\frac{4}{\gl^2}(a,y)(a,z)a
 \end{aligned}
\end{equation}
whereas
\begin{equation}
 \label{Mi2}
 (yz)^{\gt_a} = yz + \frac{2}{\gl}(a,yz)a -\frac{2}{\gl}a(yz) .
\end{equation}
Equating the two and dividing by $\frac{2}{\gl} $ gives the assertion.
\end{proof}

\begin{remark} Also recall that $a(ay) =(1-\gl)(a,y)a +\gl ay$.
For $\gl = \half$ the result reduces to \cite[Lemma 4.1, Eq. (2)]{D}.
\end{remark}

\subsection{Obtaining Frobenius $B$-axial identities from generalized Frobenius polynomials}

Another reason for studying GIs is that we also need to consider significant elements that are not idempotents, such as the elements $e,u,f$ of Lemma~\ref{tor}.
For the time being we take $C=F$, a field, although we shall work later over integral domains.

\begin{Def} A subalgebra $B$ of $A$, spanned by $\mathcal{B},$ is $1$-\textbf{hard} if every idempotent in $B$ lies in the linear span of the primitive axes of $ A$ lying in $B$.
\end{Def}

\begin{thm}\label{hard1}
 If $B$ is a 1-hard subalgebra of $A$, then any $B$-axial identity $f(E_1,\dots, E_m; X)$ which is multilinear both in the $E_i$ and the $X_j$ becomes a $B$-idempotental identity when one views the $E_i$ as idempotental indeterminates.
\end{thm}

\begin{proof} For the purposes of this proof we say $f(E_1,\dots, E_m; X)$ is a ``mixed'' idempotental identity of stage $i$ where the marked substitutions of $E_1,\dots, E_i$ are to  idempotents, and the marked substitutions of $E_{i+1},\dots, E_m$ are axial. We proceed by induction on $i.$ The case $i=0$ holds by hypothesis.
 If $a_1 $ is an idempotent, then write $a_1 = \sum \nu_{j} a'_j$ as a linear combination of axes.
 \begin{equation}
 f ( a_1,\dots, a_m ; X_1,\dots, X_n) = f \left( \sum \nu_{j} a'_j, a_2,\dots, a_m; X_1,\dots, X_n\right) = \sum \nu_{j} f ( a'_j, a_2,\dots, a_m ;X_1,\dots, X_n ) = 0,\end{equation}
and we iterate the argument, relabeling the variables at each step.
 \end{proof}

 \begin{remark}\label{dhard}
 \begin{enumerate} \item The general case for idempotental polynomials which are not multilinear can still be handled by linear algebra by taking a basis $\mathcal B$ of $B$ and expanding $f$ in terms of linear combinations of products of $b_i \in \mathcal B$ of length $k$ depending on the degree $d$ of $f$ in the $E_i.$ This would require a more cumbersome definition of $d$-hardness, which was used in our original paper.

 \item In almost all cases, hardness arises from a parametrization of the axes of Jordan type, often via \eqref{ids}.

 \end{enumerate}
 \end{remark}

\section{Axial algebras of Jordan type $\half$}

In view of Remark~\ref{MatJ}(i), we return to the important case $\gl =\half$.
We recall some Jordan algebras generated by idempotents.

 \subsection{The special Jordan algebra}

The special Jordan algebra $R^+ $ of a ring $R$ has multiplication $x\circ y = \half(xy+yx)$.

 \begin{remark}
 \label{FS} The matrix units $e_{i,i}$ are primitive idempotents of Jordan type $\half$ in $M_n(F)^+$, for any field~$F.$ The trace form $(x,y) = \tr(xy)$ is a nondegenerate normal Frobenius form on $M_n(F)^+$, since $\tr((x_1\circ x_2)x_3) = \half\tr (x_1x_2x_3)+ \half\tr (x_2x_1x_3)= \half\tr (x_1x_2x_3)+ \half\tr (x_1x_3x_2) = \tr (x_1(x_2\circ x_3))$. This restricts to a Frobenius form on any special Jordan algebra of $M_n(F)^+$.

 The special Jordan algebra satisfies the Glennie identities, which are not a consequence of the Jordan identity, cf.~\cite[Theorem~I.11.12]{J}.
 \end{remark}

Here is a special Jordan algebra of particular interest to us.

\begin{example}\label{bilf}
 We recall from \cite[\S~4]{HRS2} the Jordan Clifford algebra $Cl^J(F^{2}, \mathfrak b_\alpha)$ on the vector space $F^{2}$ with basis $v_1, v_2$ and symmetric bilinear form given by $\mathfrak b_\alpha(v_i, v_i) =2$, $\mathfrak b_\alpha(v_1, v_2) =\alpha$ for some constant $\alpha\in F.$
The Clifford algebra is spanned by $1, v_1, v_2, v_1v_2$. On the other hand, we can view the corresponding 3-dimensional Jordan Clifford algebra $Cl^J(F^{2}, \mathfrak b_\alpha)$, with basis $1, v_1, v_2$ as
  an axial algebra with generating axes
$e_1=\half (1+v_1),$ $ e_2=\half (1+v_2),$ for $\alpha\ne 2.$

More generally, the symmetric Jordan Clifford algebra $Cl^J(F^{n}, \mathfrak b)$, defined on the vector space $F^{n}$ with basis $v_1, \dots, v_{n}$ and symmetric bilinear form given by $\mathfrak b(v_i, v_i) =2$, $\mathfrak b(v_i, v_j) =\alpha_{i,j}$ for constants $\alpha_{i,j} \ne 2$ in $F,$
 has generating axes
$e_i=\half (1+v_i),$ $1\le i \le n.$
\end{example}

Note that here the idempotents are ``built into'' the construction.

\subsection{3-generated primitive axial algebras of Jordan type $\gl$}\label{3gena}\label{GSa}


 \begin{remark}
 \begin{enumerate}\eroman
 \item
Let us consider an approach to verifying that a multilinear polynomial $h(X_1,X_2,X_3,X_4)$ is an identity of a 3-generated PAJ-$\gl$ $A$.
 To check that $h(x,y,z,w)=0,$ it is enough to take $x$ to be a primitive axis~$a$ of Jordan type~$\gl$, and $y,z,w $ to be primitive axes $b,c,d$ of $A$.

Since $b^{\tau_a} = b- 2b_\gl ,$ we see at once that $\lan\lan a,b \ran \ran = \lan\lan a,b^{\tau_a} \ran \ran$ is spanned by $a,b,$ and $b^{\tau_a};$ thus we may replace $ab $ by $b^{\tau_a}.$ Likewise we may replace $ac $
 by $c^{\tau_a},$ or $bc $
 by $c^{\tau_b},$ so $y$ may also be assumed to be in $\{b,c,d\}$ and we take $y=b.$

 Next, we claim that we may replace $z$ by $c.$ This is clear by changing $c$ if $z$ is specialized to $ac,$ $bc$, or $a(bc)$
(replacing $c$ by $(c^{\tau_b})^{\tau_a}$). (Otherwise we may assume that $z=ab,$ but likewise we are done unless $w=ab,$ and then we are computing in
$ \lan\lan a,b \ran \ran,$ a Jordan algebra.)

Having proved the claim, we need merely compute $h(a,b,c,w)$, where we may assume that one of the following cases holds: \begin{enumerate}
 \item Case I: $w=a$,
 \item Case II: $w=ab$, or
 \item Case III: $w =(ab)c.$
\end{enumerate} So the assertion has been reduced to three computations, which could be checked via \cite[Table~I]{GSS} or \cite{DRS}.

\item (A heuristic computer-free approach to the Gorshkov-Staroletov theorem that 3-generated PAJs are Jordan algebras)
 Here is a quick application. In the notation of (i), the proof of \cite[Proposition~4, p. 17]{GSt}, that each 3-generated PAJ-$\half$ is a Jordan algebra,
 could be finished off fairly quickly with~\cite{DRS}. Full computations would still add a few pages to the paper, so we content ourselves with a heuristic argument. We consider the case
 of a unital algebra with nondegenerate Frobenius form. First of all, in view of Fact~\ref{br}(xiii), in characteristic not dividing 30, one need only consider the multilinearization of 4-power-associativity \eqref{PA1} given in \cite[Eq (3)]{Ko}, namely

\begin{equation}\label{powa1}
 \begin{aligned}
 h(x,y,z&,w):= - 4[(xy)(zw) + (xz)(yw) + (xw)(yz)]
\\& + x[y(zw) + z(wy) + w(yz)] + y[x(zw) + z(wx) + w(xz)]
\\& + z[x(yw) + y(wx) + w(xy)] + w[x(yz) + y(zx) + z(xy)]
 \end{aligned}
\end{equation}
for $x,y,z,w \in A $.

Specifically, to check that any 3-generated PAJ-$\half$ satisfies~\eqref{powa1}, perhaps the best way to compute is via the left multiplication maps $L_y(z) = yz.$ Case I is easy. The other two cases are too lengthy to write out here, but are computable directly.

In Case II, $w = ab,$
\begin{equation}\label{powa2}
 \begin{aligned}
 h&(a,b,c,ab):= - 4[(ab)(c(ab)) + (ac)(bab) + (a(ab))(bc)]
\\& + a[b(c(ab)) + c((ab) b) + (ab) (bc)] + b[a(c(ab) ) + c((ab) a) + (ab) (ac)]
\\& + c[a(b(ab) ) + b((ab) a) + (ab) (ab)] + (ab) [a(bc) + b(ca) + c(ab)].
 \end{aligned}
\end{equation}
can be rewritten
\begin{equation}\label{powa2a}
 \begin{aligned}
 \tilde h(L_a,L_b,&L_{ab}) := - 4(L_{ab}L_{ab} + L_{bab} L_a + L_{a(ab)}L_b)
\\& + L_a(L_b L_{ab} + L_{(ab) b} + L_{ab} L_b) + L_b(L_aL_{ab} + L_{ab a} + L_{ab} L_a)
\\& + L_{a(bab) } +L_ {b(ab a)} + L_{(ab)(ab)} + L_{ab}(L_aL_b + L_bL_a + L_{ab}).
 \end{aligned}
\end{equation}

The point is that there are formulas for $L_{aba}$ and $L_{bab}$ in \cite[Proposition~1.2 and Proposition~2.1]{DRS}, which would yield \eqref{powa2a}.

Case III: $w = (ab)c.$
Since the terms are linear in $w$, in view of Cases I and II, we could replace $w$ by a linear combination of $ac, bc,$ and $b^{\gt_a}c.$ The first two are covered by Case II, so we are left with $w = b^{\gt_a}c$. Switching $b$ and $b^{\gt_a}c$ leaves us with $h(a,b^{\gt_a},c,bc),$ which becomes $h(a,ab,c,bc)$. This expression is somewhat trickier than before, but manageable.

In summary,  this does not constitute a full independent proof of the Gorshkov–Staroletov theorem, since the remaining computations are not supplied.
\end{enumerate}
 \end{remark}

When $x,y$ in \eqref{powa1} are specialized to the same idempotent, we obtain

\begin{prop}
 Every PAJ-$\half$ satisfies the $A$-axial identity
 \begin{equation}
 \begin{aligned}
 4[E_1& (X_1X_2) + (E_1X_1)(E_1X_2) + (E_1X_2)(E_1X_1)] =
\\&2E_1[E_1(X_1X_2) + X_1(E_1X_2) + X_2(E_1X_1)]
\\& + X_1[2E_1(E_1X_2) + E_1X_2] + [2E_1(E_1X_1) + E_1X_1] X_2
 \end{aligned}
 \end{equation}
\end{prop}
\begin{proof}
 To verify the identity, in view of Proposition~\ref{Miy6}, we may assume that $A$ is 3-generated, but then $A$ is Jordan by \cite[Proposition~4, p. 17]{GSt}.
\end{proof}

\subsection{Solid subalgebras: The theorems of  \cite{D} and \cite{GSS}}

Having the framework in place, we want to interpret the following theorems in \cite{D} and \cite{GSS} about solid subalgebras of a PAJ-$\half$, studied in terms of its 2-generated subalgebras. We are not presenting full proofs, but only indicating how the ideas might be understood in terms of axial identities.

\begin{Def}(\cite{D})
 A subalgebra of a PAJ-$\half$ is \textbf{solid} if all of its primitive idempotents $\ne 0,1$ become axes when extending the base field to an algebraically closed field.
\end{Def}

 \cite[Definition~5.1]{GSS} also requires that the primitive axes be of Jordan type $\half.$ We reconcile the discrepancy below.
We want to address the following major results about solid subalgebras of PAJ-$\half$ over a field $F$ of order $>3$.

\begin{thm}\label{GSSD} Let $A$ be a PAJ-$\half$, with a normal Frobenius form $(\phantom{w},\phantom{w})$.
 \begin{enumerate}\eroman
 \item For primitive axes $a, b$ of Jordan type $\half$, the
subalgebra $\lan\lan a,b \ran\ran$ is solid whenever
\begin{enumerate}
 \item (\cite[Theorem 1]{GSS}) $(a, b) \notin \{0, \frac{1}{4},1\},$
 \item (\cite[Theorem 2]{GSS}) $\charc\ F =0$ and $(a,b) \ne \frac{1}{4}$,
 \item (\cite[1.1]{D}) $\lan\lan a,b \ran\ran$ is 2-dimensional.
\end{enumerate}
 \item (\cite[1.3]{D}) If $\charc\ F \ne 3$ and
$\lan\lan a,b \ran\ran$ is solid for all $a,b\in \mathcal{E}$ then $A$ is a Jordan algebra.
 \end{enumerate}
\end{thm}

\begin{remark}
 To tell whether a subalgebra is solid,
we need to check the axial Frobenius identities of Remark~\ref{summ}, which say that $x$ is an axis of Jordan type. The method espoused here handles this verification directly.
\end{remark}

These notions can all be reformulated in terms of idempotental identities and axial
identities.

\begin{remark}\label{ide} Let $A$ be a PAJ-$\half$, with a sub-PAJ $\lan\lan a,b \ran\ran$.
 \begin{enumerate}\eroman
 \item $\lan\lan a,b \ran\ran$ is solid in the sense of \cite{D} if and only if $A$ satisfies the primitive $\lan\lan a,b \ran\ran$-idempotental identities from Example~\ref{prim0} and Lemma~\ref{axfree}.

 \item $\lan\lan a,b \ran\ran$ is solid in the sense of \cite{GSS} if and only if $A$ satisfies the primitive $\lan\lan a,b \ran\ran$-idempotental identities implicit in \S\ref{axid}.

 \end{enumerate}
\end{remark}

Thus, Theorem~\ref{GSSD}(i) is about proving that certain primitive $V$-idempotental identities are $V$-axial identities, and Theorem~\ref{GSSD}(ii) is about proving that certain $\mathcal{E}$-axial identities imply the Jordan identity. This was one of the key steps in \cite{GSS,D}, as seen in:

\begin{Facts}\label{inf2} $B$ is 1-hard for $B = \lan\lan a,b \ran\ran\subseteq A$ in the following situations:
 \begin{enumerate}
 \item It is easy to see that for $F$ infinite, infinitely many values of $\gve$ in \eqref{ids} yield primitive axes of $A$, which occurs when
 \begin{enumerate}
 \item \cite[Lemma~3.9]{GSS} $A$ is the Jordan algebra of a Clifford algebra.
 \item \cite[Corollary~5.4]{GSS} $(a,b)\notin \{0,1,\frac{1}{4}\}$.
 \item \cite{GSS,D} $B$ is flat.
 \end{enumerate}
 \item $F$ has characteristic $0$ (as noted in the proof of \cite[Corollary~6.5]{GSS}), since $\gt_a \gt_b$ has infinite order. \end{enumerate}
\end{Facts}

\subsubsection{The key step in proving the theorem of Desmet}

Desmet~\cite{D} showed
 that a primitive axial algebra of Jordan type $\tfrac{1}{2}$ is a Jordan algebra if and only if every 2-generated, subalgebra is solid.

First of all, \cite[Lemma~4.1]{D} can be stated in terms of axial identities.
 We give an axial-identity approach, which is designed to provide the key implication in the toric case over an infinite field.

 In principle, we follow the approach of \cite{D}.\footnote{This approach is intended to replace \cite[Corollary~5.5]{D}.} 
      Since $A$ is spanned by axes, it is enough, in view of Facts~\ref{br}(ix),  to prove that $[L_a,L_b]$ is a derivation for any axes $a,b$, which is easily translated into a multilinear axial identity
      \begin{equation}\label{deriv}
          \begin{aligned}
              g(E_1,E_2; X_1,X_2) =& E_1(E_2(X_1X_2))- X_1 (E_1(E_2X_2)) - (E_1(E_2X_1))X_2\\& -  E_2(E_1(X_1X_2))+ X_1 (E_2(E_1X_2))+(E_2(E_1X_1))X_2
          \end{aligned}
      \end{equation}

    \begin{conj}
        Suppose $B = \langle c,d \rangle$ is a toric pair in an algebra spanned by axes, and $g(E_1,E_2;X_1,X_2)$ is an idempotental multilinear polynomial which is symmetric in   $X_1,X_2,$ such that $g(a_1,a_2,c,c)=0$ for all axes $a_1,a_2$ and all $c\in B.$ Then is $g$ is an axial identity of $A$?
    \end{conj}
    This would be a useful result, so we describe an approach towards a proof. It suffices to specialize $X_1,X_2$ to toric axes.  Since $\langle c,d \rangle$ is toric, we know by hypothesis that all nontrivial idempotents of~$\lan\lan c,d \ran \ran$ have the form $$  x_\gve  = \gve  e  + \gve^{-1} f  + \half u , \quad \gve \in F^\times,$$ where infinitely many are axes.

    Write  $\gve_1\gve_2 g\left( \bar E_1, \bar E_2,\gve_1^2   e +      f + \half \gve_1 u,   \gve_2^2 e +      f + \half \gve_2  u\right) = \sum h_k \sum _{i+j=k} h_{i,j}$, where $h_{i,j}:= h_{i,j}( \bar E_1, \bar E_2, u,e,f).  $
    We want to verify that each $h_{i,j}=0.$

    But taking $\gve_1=\gve_2= \gve$ yields $\gve^2 g\left( \bar E_1, \bar E_2,\gve^2   e +      f + \half \gve u,   \gve^2 e +      f + \half \gve  u\right) = \sum_{k=0}^4 \sum _{i+j=k}  \gve ^k\sum h_{i,j} . $
    Since this holds
    for infinitely many $\gve,$ we have $\sum _{i+j=k} \sum h_{i,j} =0$ for all $k.$ Now $h_{i,j}=h_{j,i} $ since $g$ is symmetric. Hence we need to show each $h_{i,j}= \half(h_{i,j}+h_{j,i} )=0.$
    For $k = 0$, we have $h_{0,0} =0.$   For $k = 1$, we have $2h_{1,0}  =0,$ so   $h_{1,0}=0.$ For $k = 2$, we have $ 2h_{2,0} +h_{1,1} =0.$
     For $k = 3$, we have
$h_{2,1} =0$.
    For $k = 4$, we have
$h_{2,2}=0$.
We want each $h_{i,j} = 0$ to conclude the proof, so it remains to show that $h_{1,1} = 0 $.

We do not see yet how to do this, we but note that since $g$ is linear in $E_1,$ we can replace $\bar E_1$ by any element in $\langle \bar E_1, \bar E_2,  c\rangle,$ for example a word in $c  \bar E_1.$  Since any 3-generated PAJ has finite dimension, this might conclude the argument.

Note that
 $g$ of \eqref{deriv} would be an axial identity for Jordan algebras by Osborne's theorem, since $a_1,a_2,c$ generate a Jordan algebra, by   \cite{GSt}. Hence $g$ would be an axial identity for $A$

 \subsubsection{Proposed approach to the theorem of Gorshkov, Shpectorov, and Staroletov}

The following argument is conditional on extending 1-hardness of $B$,  by Facts~\ref{inf2}, to the required $d$-hardness
indicated in Remark~\ref{dhard}. If so, each primitive idempotent would satisfy the $B$-axial identity of Jordan type $\half,$ making it an axis.

\subsubsection{Fusion rules on hard subalgebras}

Knowing that idempotents in hard subalgebras are primitive axes, the next question is whether they are of Jordan type.

\cite[Proposition 6.3]{GSS} shows in characteristic 0 that the anisotropic primitive axes are of Jordan type. Theorem~\ref{hard1} shows this in general, for hard axial subalgebras. 

\section{Generic axial algebras}\label{gaa}

 Besides providing a framework for structural results, idempotental polynomials yield a rich source of examples of axial algebras. We have the following hierarchy of constructions:

 $$
\text{relatively free}
\longrightarrow
\text{finite subdirect decomposition}
\longrightarrow
M\text{-primitivity}
\longrightarrow
\text{generic PAJs}.
 $$

We start with the ``relatively free'' marked axial algebra. Since there may be torsion over the coefficient ring, we take a subdirect decomposition into finitely many generically local components.
We impose $M$-primitivity by adjoining scalar parameters for the 1-components. These components provide the desired generic MPAJ-$\gl$ algebras upon passage to quotient fields.

\begin{remark}
 If $U$ is torsion-free over some integral domain $\overline{C},$
 and $K= Q(\overline{C})$,
then $U^{K}$ is an axial algebra over a field. This is why we would like torsion-free generic axial algebras. \end{remark}

\cite{HRS1} defined ``universal PAJs'' using a Frobenius form, but here we take a varietal approach in which the Frobenius form is not needed initially, thereby enabling us to treat Shpectorov's Question \ref{Spe} conceptually. To carry this out, we extend our base ring $C$ to contain
commuting indeterminates, then mod out by an ideal of axial identities which yields a specialization property.
\subsection{Relatively free axial algebras}

\begin{Def}\label{fr1}
 \begin{enumerate}
 \item The \textbf{free nonassociative algebra} $C\{Y\}$ in $N$ generators is generated by ``free'' indeterminates $Y_i,\ 1\le i\le N$.

 \item A \textbf{relatively free algebra} is an algebra $C\{Y\}/\mathcal I$ where $\mathcal I$ is a T-ideal.

 \item The \textbf{free $\gl$-axial algebra} $C\{E\}$ in $N$ idempotents   over $C$ is obtained by modding out the relations $\mcI = \{E_i^2 -E_i: 1\le i\le N\}$ and the idempotental identities $E_i(E_i(E_ig) - (\gl +1)E_ig + \gl g)$, for all $1\le i\le N$, $g\in C\{E\}$. Thus, writing $E_i$ for the image of $Y_i$, we have the free nonassociative algebra generated by ``free'' idempotents $E_i^2 = E_i,\ i\in I$\footnote{Our original definition also had indeterminates $X_j$, but this raises the complication that the elements $X$ are not generated by axes.}

 \item An \textbf{axial T-ideal} $\mathcal I$ of $C\{E\}$ is an ideal invariant under all marked substitutions of the $E_i$ in $C\{E\}$ to axes, in other words invariant under all marked endomorphisms fixing $C$.
 \item The {axial T-ideal} $\langle \mathcal S \rangle$ \textbf{generated by} a subset $\mathcal S$ of $C\{E\}$ is the intersection of all axial T-ideals containing $\mathcal S$.
 \item A \textbf{relatively free axial algebra} is an algebra $C\{E\}/\mathcal I$ where $\mathcal I$ is an axial T-ideal.
 \end{enumerate}
\end{Def}

\begin{remark}\label{Tfacts}
 As in the usual PI-theory,
 \begin{enumerate}
 \item \cite[Theorem~2.1.8]{Row}, when $\mathcal I$ is an axial T-ideal of $C\{E\}$, $\mathcal I$ corresponds to the pure axial identities of $C\{E\}/\mathcal I$.

 \item If every element of $\mathcal S$ is an axial identity of a marked algebra $A$, then so is every element of $\langle \mathcal S \rangle$.

 \item For $c\in C,$ $\langle c \mathcal S \rangle = c \langle \mathcal S \rangle$.
 \end{enumerate}
\end{remark}


\begin{example}
\begin{enumerate}\eroman
 \item
The relatively free 3-generated Jordan algebra $U_{\operatorname{Jor}}$ over generators $Y_1,Y_2,Y_3$ is obtained when we mod out the T-ideal generated by the Jordan identity \eqref{nonc}.

\item
The relatively free axial Jordan algebra $U_{\operatorname{Jor}}$ over axes $E_1,E_2,E_3$ is obtained when we mod out the axial T-ideal generated by evaluations of $((xx)y)x - (xx)(yx)$, i.e., yielding the Jordan identity \eqref{nonc}.

\item Other relatively free axial Jordan algebras are obtained by modding out extra identities, including the Glennie identities.

\end{enumerate}
\end{example}

When $\overline{C}$ is an integral domain,
 as before, write $ K:=Q(\overline{C})$, and $U^{K}$ for $U\otimes_{\overline{C}} {K}.$

 \begin{remark}\label{bas1} Pick a spanning set $\mathcal{B}$ of $U^{K}$ over ${K},$ containing $\mcE = \{E_1,\dots,E_N\}.$
 \begin{enumerate}
 \item We can choose $\mathcal{B}\subseteq U $, and so any element of $U^{K}$ is a linear combination of the spanning set with coefficients in $K.$ Suppose $S = \{b_i: i\in I\}$ is a finite subset of $\mathcal{B}$. Then writing $b_i b_j = \sum c_{i,j,k}b_k,$ we can clear denominators and replace the elements of $S$ by scalar multiples such that $c_{i,j,k}\in \overline{C}.$
 \end{enumerate}
 \end{remark}

 \subsection{Marked generic axial algebra}

\begin{Def}
 \begin{enumerate}
 \item An algebra $U$ generated by marked elements $\{E_i:i\in I\}$ is \textbf{generic} in a marked variety $\mathcal V$ if for every algebra $A$ in $\mathcal V$ and marked elements $\{a_i:i\in I\},$ there is a  marked homomorphism $\psi: U \to A$ sending $E_i\to a_i $, $\forall i \in I.$ 

 \item A set of algebras $U_1,\dots, U_k$ with each $U_\ell$ generated by marked elements $\{E_{i,\ell} :i\in I\}$ is \textbf{generic} in a marked variety $\mathcal V$ if for every algebra $A$ in $\mathcal V$ and marked elements $\{a_i:i\in I\},$ there is some $U_\ell,$ $1\le \ell \le k,$ 
 and a  marked homomorphism $\psi: U_\ell \to A$ sending $E_{i,\ell}\to a_i $, $\forall i \in I.$ 
 \item More specifically, if $\mathcal V$ is an axial variety and $U_\ell,$ $1\le \ell \le k$ are generically local axial algebras over integral domains, we call $\{U_\ell: 1\le \ell \le k\}$ \emph{the generically local} axial algebras of $\mathcal V$.
 \end{enumerate}
\end{Def}

In what follows, most of our generic algebras will be constructed as relatively free axial algebras. Here is an exception, where we expand the language to contain a bilinear form.
 \begin{example}
 The generic Jordan Clifford algebra over a vector space $V$ with basis $\{b_1,\dots, b_t\}$ is obtained by applying the generic trace bilinear form $(b_i,b_j)= \chi_{i,j}$ for commuting indeterminates $\chi_{i,j}, 1 \le i \le t$. One nice feature of this construction is that one obtains specializations merely by specializing the $\chi_{i,j}$.
 \end{example}

\subsection{Constructing marked generic $\gl$-axial algebras}

We start with the following construction (which works more generally for all marked varieties). We still assume that $\gl(1-\gl)$ is invertible in $C.$

 For $M> N$ to be specified, we index all words in $E_1,\dots, E_N$ of length $\le M$ as $Y_1,\dots, Y_{\theta(N,M)}\}$, where $Y_i = E_i$ for $1\le i \le N,$
and expand $C$ to the polynomial ring $C^\natural = C[\Xi]:=C[\xi_{i,j} \ : \ 1\le i \le N, \ 1\le j\le \theta(N,M)]/\langle \xi_{i,i}-1, \ 1\le i\le N \rangle$ in the commuting indeterminates $\xi_{i,j}$. 

\begin{remark} $C^\natural$ is Noetherian, by the Hilbert basis theorem, so the tools of affine algebraic geometry are readily available. (This is the reason we introduce~$M,$ taking it large but finite.)

 We could formulate the theory for arbitrary $S$-axes, but take $S = \{\gl\}$ and work with PAJ-$\gl$ to simplify the discussion.
The situation is clearest when $N$ is $\gl$-axially admissible, i.e., if each $N$-generated PAJ-$\gl$ satisfies $A = A^{(M)}$, which is seen without undue difficulty for $N\le 4,$ by \cite{DRS}.
For $N\ge 5$ the situation is more delicate.
\end{remark}

 Our technique will be to quotient $C^\natural\{E\}$ by an ideal of axial identities. This is a $C^\natural$-algebra but need not be torsion-free, an obstacle which we need to cope with.

\begin{Def}
 \begin{enumerate}
 \item List the words $Y_j, 1\le j \le \theta(N,M)$ as before, and $U_{\gl,N,M}$ is the relatively free algebra $C^\natural\{E\}$ modulo the axial T-ideal generated by $E_i(E_i(E_iY_j))- (\gl +1) E_i(E_iY_j)+\gl E_iY_j,$ and $ \gl E_iY_j -E_i(E_iY_j) - ({\gl-1})\xi_{i,j}E_{i}$, for all $1\le i \le N$ and all $1\le j \le \theta(N,M),$ cf.~\eqref{ax10b}.

 \item $U_{\gl,\pp}$ is $U_{\gl,N,M}$ modulo the axial identities common to all PAJ-$\gl$ over fields which contain a homomorphic image of $C$.

 \end{enumerate}
\end{Def}

\begin{lemma}\label{lem4.6} Let $U = U_{\gl,\pp}.$ Suppose $a\in \mcE.$ $U_1(a) $ is a commutative, associative algebra with unit element $a,$ whose nilpotent radical is 0.
\end{lemma}
\begin{proof}
 In any homomorphic image PAJ-$\gl$ $A$ over a field $F$, $ {A_1(\bar a)} $ is commutative, associative, and the nilpotent radical is 0 (since it is $F \bar a$). Hence $U_1(a)$ is commutative and associative, since these properties are expressed in terms of polynomial identities. Furthermore, if $f$ is in the nilpotent radical of $U_1(a)$ then $f^k = 0$ for some~$k,$
 implying $\bar f ^k = 0$ in any PAJ-$\gl$ $\bar A$, so $\bar f = 0$, i.e., $f=0$ holds in $U$.
\end{proof}
\begin{remark}

This does not say even for $N=4$ that $U_{\gl,\pp}$ is a PAJ. In fact $U_{\gl,\pp}$ may fail to be torsion-free over $C^\natural.$ $U_{\gl,\pp}$ need not have a Frobenius form, which has been our tool to determine when an axis is primitive. Thus we would want to find some way using generalized identities to prove primitivity. This seems highly unlikely.
\end{remark}

\subsection{Removing torsion}
 Here is a way to produce finitely many generic axial algebras from a relatively free axial algebra $U= C^\natural\{E\}/\mcI$, denoted $\overline{C^\natural\{E\}}$, viewed over $C^\natural$.

We view the elements of $U$ as pure polynomials in the $\bar E_i.$ The \textbf{degree} of $\bar y =f(\bar E_1,\dots,\bar E_N)\in U$ is the lowest degree of $f$ which can be used to present $\bar y.$

\begin{lemma}
 If $\mcI$ is a T-ideal of $U$, then $\Ann_{C^\natural} \mcI$ is an ideal of $C^\natural,$ so $\Ann_U (\Ann_{C^\natural} \mcI)$ is a T-ideal of $U.$
 \begin{proof}
 Just push through the substitutions.
 \end{proof}
\end{lemma}

Recall $\mcT$ from Definition~\ref{annid}, replacing $A$ by $U$ and $C$ by $C^\natural$.
 \begin{Def}
 \begin{enumerate}
 \item A \textbf{torsion polynomial} of $U$ is an element $f\in U$ (i.e., a pure idempotental polynomial) for which $\Ann_{C^\natural}(f)$ contains a regular element of $C^\natural$.

 \item An ideal $\mcI$ in $\mcT$ is \textbf{torsion} if $\Ann_{C^\natural}\mcI $ contains a regular element.

 \item $\mcT-\tor (U)$ is the sum of torsion ideals in $\mcT$.

 \item $d-\mcT-\tor (U)$ is the submodule of $U$ consisting of those pure axial polynomials having degree $\le d$ in the $\bar E_i$ which are torsion.
 \end{enumerate}
 \end{Def}
\begin{lemma}
 $\mcT-\tor(U/\mcT-\tor(U)) =0.$
 \end{lemma}
 \begin{proof} If $h_i f_i=0$ for regular $h_i \in C^\natural$, $i = 1,2,$ then $h_1h_2(f_1+f_2)=0$ and $h_i(g f_i)= 0 = h_i(f_ig)$.

 If $h_1(\bar f +\mcT-\tor (U))\in \mcT-\tor (U),$ then $h_2h_1 \bar f$ is an axial identity, so $\bar f\in \mcT-\tor (U).$
 \end{proof}

 \begin{Def}
 $U$ is \textbf{generically $\mcT$-local} if $C^\natural$ contains a unique maximal $\mcT$-annihilator~$P$ of $U$.
 \end{Def}

\begin{thm}\label{TFa}
 Suppose $C$ is Noetherian, and let $\mcI$ be an $\mcT$-ideal. Take $U = C^\natural\{E\}/\mcI. $
 \begin{enumerate}
 \item $U$ is a subdirect product of finitely many generically $\mcT$-local axial algebras $U_\ell$.

 \item Notation as in Lemma~\ref{lem4.6}, for any $a\in \mcE,$ define $U_1^{(M)}(a) $ to be the span of $\{y_1 : \operatorname{length}(y) \le M $\}. Then $U_1^{(M)}(a) $ is a finite module over~$C^\natural.$
 \end{enumerate}
\end{thm}
\begin{proof}
 We prove the theorem by Noetherian induction on the $\mcI$-annihilator ideals, since $C^\natural$ is Noetherian.    Namely,
we are done unless some ideal in $\mcT$ has an annihilator in $C^\natural$. Take a
 maximal $\mcT$-annihilator $P_1$ in $C^\natural$.
 Let $ \mcI_1= \Ann_U (P_1).$ Thus $P_1 = \Ann_{C^\natural} (\mcI_1)$.

We are done unless there is a second maximal $\mcT$-annihilator $P_2\ne P_1$, say of $\mcI_{2}$. Then the canonical homomorphism $\Psi: U \to U/\mcI_{1}\times U/\mcI_{2}$ is an injection since any element of the kernel $\mcI_{1}\cap \mcI_{2}$ would be annihilated by $P_1+P_2,$ an ideal of $C^\natural,$ contrary to the maximality of $P_1.$

 For $\ell = 1,2$, each of the $ U/\mcI_{\ell} $ is an  axial algebra over    $C^\natural$,  so by Noetherian induction on the $P_i \supset \Ann_{C^\natural}U$, is a finite subdirect product of generically $\mcT$-local  axial algebras which are faithful over suitable prime images of $C^\natural,$ and we are done.

 (ii) In view of Lemma~\ref{axfree}, $U_1^{(M)}(a)$ is isomorphic to a submodule of $\prod _{1\le \ell \le k} U_\ell^{(M)}$, which is finite as a $C^\natural$-module.
\end{proof}

This leaves us with the question:

\begin{ques}
 \begin{enumerate}
 \item If $U$ is generically $\mcT$-local, is $U$ then $\mcT$-torsion-free over $ C^\natural/\Ann_{C^\natural}(U)$?

 \item More specifically, if $h(\Xi)\mcI \subseteq \mcI_1$ for $\mcI, \mcI_1 \in \mcT$ and $h\notin \Ann_{C^\natural}(U),$ then is $\mcI = 0$?
 \end{enumerate}
\end{ques}

We have a partial answer over infinite fields. We identify $f(E;X)$ as an axial identity when $f(E;g)$ specializes to 0 for all $g\in U.$

 \begin{ques}
 If $C$ is an infinite field $F$, and $h(\Xi)f(E_1,\dots, E_m)$ is a pure axial identity of a finite-dimensional PAJ $A$ over $F$. Then must either $hA$ be $0$ or $f$ be an axial identity of $A$?
 \end{ques}

 Equivalently, taking a generically local algebra $U$
 lying over $A,$ if $hf = 0$ in $U,$ then is $h= 0$ or $f=0?$ We would like some sort of density argument on the marked specializations.

\subsection{Generic Matsuo algebras}

 We turn to Matsuo algebras, which are easier to study since primitivity is built in by Lemma~\ref{M1}.

 \begin{Def} $M_\gl$ is $U_{\gl,\pp }/{\mcI}_{\operatorname{Mat}}$, where ${\mcI}_{\operatorname{Mat}}$ is the T-ideal of pure idempotental polynomials arising from \eqref{Mat2a} and \eqref{Mat3a}.
\end{Def}

\begin{thm}\label{gen2} There is a finite set of generically $(\mathcal T)$-local relatively free marked algebras in the variety of 4-generated Matsuo algebras of Jordan type $\gl$, i.e., for which every 4-generated Matsuo algebra is isomorphic to a homomorphic image of one of them.
\end{thm}
\begin{proof} $\mathcal I_{\operatorname{Mat}}$ is the set of pure axial identities of the marked variety of Matsuo algebras, by Lemma~\ref{Matl}, so
 Theorem~\ref{TFa} (taking $\mcT$ to be the axial T-ideals) is applicable.
\end{proof}

 The assertion also can be proved for $N$-generated Matsuo algebras of Jordan type $\half$ for $N>4$, but relying as of now on \cite{CuH}. The proof does not state how many algebras are needed, but one can obtain a bound on the number of possible normal Frobenius forms, since the values are in~$\{ 0, 1, \frac{\gl}{2}\}.$

 \subsection{Making axes $M$-primitive}

 We fix $M$ sufficiently large. (This refinement is automatic for $N= 4.$) Under marked specializations, $\xi_{i,j}$ is to be specialized to the scalar $\varphi_{a_i}(y
_j)$ of Definition~\ref{not1}.

 \subsubsection{Marked generic MPAJ-$\gl$ algebras}

\begin{remark}\label{Fr2}
 \begin{enumerate}\eroman
 \item The base ring $C^\natural$ need not be an integral domain, so the situation is delicate.

 \item The generic axes are $M$-primitive axes whose left multiplication maps satisfy a polynomial of degree~3, since one specializes down to the axes generating a PAJ-$\gl$.

 \end{enumerate}
\end{remark}

 In Theorem~\ref{TFa}, taking $U=U_{\gl,\pp},$ let $U_{\gl,N_\ell;M-\operatorname{prim}_\ell},$ $1 \le \ell \le k$, denote the corresponding $N_\ell$-generated axial algebras $U_\ell$. (Note that some $E_i$ could be mapped to 0, so $N_\ell \le N.$)

\begin{thm}\label{gen1} Each $U_{\gl,N_\ell;M-\operatorname{prim}_\ell}$
is of Jordan type~$\gl$, and is generically $\mcT$-local over $C_\ell = C^\natural/P_\ell$.

 There is a canonical homomorphism $\Psi: U_{\gl,\pp}\rightarrow \prod _{\ell =1}^k U_{\gl,N_\ell;M-\operatorname{prim}_\ell}$, whose kernel is $\mcT-\tor(U_{\gl,\pp})$.
\end{thm}
\begin{proof}
 By Remark~\ref{quo}(i), we use the technique of Theorem~\ref{TFa} (taking $\mcT$ to be the collection of T-ideals) to provide the relatively free axial algebras $U_{\gl,N_\ell;M-\operatorname{prim}_\ell}$, which are generated by the nonzero images of the $\bar E_i,$ $1\le i\le N,$ which are axes of Jordan type~$\gl$, and having the universal property that any MPAJ-$\gl$ over $F$ is a marked homomorphic image of one of them.
 We conclude with Remark~\ref{quo}.
\end{proof}

Note: For $\gl\ne \half,$ only homomorphic images of Matsuo algebras can arise, and these have already been handled.

Next, we want to decompose further, and arrive at genuine PAJs over a field. In the next result, in the set-up of Theorem~\ref{gen1} we could take $U$ to be one of the $U_{\gl,N_\ell;M-\operatorname{prim}_\ell}$, and write $U_{\gl,N_\ell;M-\operatorname{prim}_\ell}$ for $U_\ell$.

\begin{thm}\label{gen11}
 In Theorem~\ref{TFa} take $\mcT$ to be the collection of all proper ideals of $U$.
 \begin{enumerate}\eroman
 \item $C_\ell$ has a unique maximal $\mathcal T$-annihilator associated with $U_\ell.$
 \item Write $K_\ell$ for the field of fractions of the integral domain $C_\ell$, and write $\pi_\ell$ for the composition of the projection from $U $ to the $\ell$ component, $U_\ell$, followed by the inclusion $U_{\gl,N_\ell;M-\operatorname{prim}_\ell} \hookrightarrow U_{\gl,N_\ell;M-\operatorname{prim}_\ell}^{K_\ell}$.
 If $c=\pi_\ell (E_i)\ne 0$ for $E_i\in \mcE,$ then $c$ is an $M$-primitive axis in $U_{\gl,N_\ell;M-\operatorname{prim}_\ell}^{K_\ell}.$

 \item $U_{\gl,N_\ell;M-\operatorname{prim}_\ell}^{K_\ell}$ is an MPAJ-$\gl$ over the field $K_\ell,$ for each $\ell.$

\item Each axial identity of $U_{\gl,N_\ell;M-\operatorname{prim}_\ell}$ is an axial identity of
 $U_{\gl,N_\ell;M-\operatorname{prim}_\ell}^{K_\ell}$.
 \item The algebras $U_{\gl,N_\ell;M-\operatorname{prim}_\ell}$, $1\le \ell \le k,$ form a set of at most $k$ generic  algebras for MPAJ-$\gl$ over~$C,$ in the sense that each MPAJ-$\gl$ over $F$ is a specialization of one of these. (However, the individual $U_{\gl,N_\ell;M-\operatorname{prim}_\ell}$ need not be relatively free.)

 \item Now assume $N$ is $\gl$-axially admissible, and write $U_{\gl,N,\operatorname{prim} }$ for $ U_{\gl,N,M-\operatorname{prim} }$ where $M = \Phi(N)$. Then
 $\prod _{\ell =1}^k U_{\gl,N_\ell;M-\operatorname{prim}_\ell}^{K_\ell}$ has a normal Frobenius form.
 \end{enumerate}
\end{thm}
\begin{proof}
(i) Generically $\mcT$-local translates to this condition.

(ii) When $c\ne 0,$ it satisfies the idempotental polynomials
$E_i(E_i(E_iY_j))- (\gl +1) E_i(E_iY_j)+\gl E_iY_j,$ and $ \gl E_iY_j -E_i(E_iY_j) - ({\gl-1})\xi_{i,j}E_{i}$, for all $1\le i \le N$ and all $1\le j \le \theta(N,M),$ cf.~\eqref{ax10b}
describing $M$-primitivity, so the axis $c$ is $M$-primitive.

(iii) By (ii) applied to each $E_i$ in $\mathcal{E},$ for which $\pi_\ell(E_i)\ne 0,$ $U_{\gl,N_\ell;M-\operatorname{prim}_\ell}^{K_\ell}$ is an MPAJ-$\gl$ over the field $K_\ell$.
Any axial specialization from $U\to A$ passes through some $P_\ell$, inducing $U_{\gl,N_\ell;M-\operatorname{prim}_\ell}\to A$.

(iv) If ${C^\natural}/P_\ell$ is a finite ring then it is already a field. Otherwise, the assertion is proved in \cite[Corollary~2.3.32]{Row}.

(v) Follows at once from (iv).

(vi) Since $M=\Phi(N)$, we have $A=U_\ell^{K_\ell}$; hence $M$-primitivity is ordinary primitivity. Each $U_{\gl,N_\ell;M-\operatorname{prim}_\ell}^{K_\ell}$ is endowed with the normal Frobenius form by \cite{HSS1}, since it is a PAJ-$\gl$; we piece together the components.

\end{proof}

 \subsubsection{Generic nonsingular MPAJ-$\gl$ algebras}

We can talk of ``nonsingular'' even without a Frobenius form, in view of Remark~\ref{rady}(iii). We follow a Baer-type construction to build up elements inside the radical.

   \begin{Def}\label{gen7}
  \begin{enumerate}
     \item  The \textbf{ nonsingular Baer radical} $\mathcal{R}(A)$  of an $\gl$-axial algebra $A$ spanned by the axes $\tmcE$ is defined via  transfinite induction as follows:
$\mathcal{R}_0$ is the ideal generated  by all $y\in A$ such that $a(a y) = \gl a y$ for all   $a\in {\tmcE}$.
Inductively, for any ordinal $i,$ $\mathcal{R}_{i+1}\supseteq \mathcal{R}_i$ is the ideal generated  by all $y\in A$ such that $a (ay) - \gl (ay) \in \mathcal{R}_i$ for all  $a\in {\tmcE}$.  For any limit ordinal  $\delta,$ $\mathcal{R}_{\delta} = \bigcup_{i<\delta} \mathcal{R}_{i}$.  $\mathcal{R}(A) := \bigcup _i \mathcal{R}_i;$ the chain must stabilize  before $|A|^+$, because each step increases the cardinality.

\item   Write $\bar A = A/\mathcal{R}(A)$.
\end{enumerate}
\end{Def}

\begin{lemma} \begin{enumerate}\eroman
    \item In the above notation
    $ \bar A$ is nonsingular.
    \item $\overline{U_{\gl,\pp}}$ has an axial specialization to every nonsingular $M$-primitive axial algebra.
\end{enumerate}
\end{lemma}
\begin{proof}(i) Suppose   $a(ay) -\gl ay\in \mathcal{R}_{i_a}.$
Letting $i = \sup \{i_a :a\in\tmcE\},$ we see  that
 $a(ay) - \gl ay \in  \mathcal{R}_i \subseteq \mathcal{R}.$ Hence     $y \in \mathcal{R}_{i+1},$
 implying $\bar y = 0.$

  (ii)  Every  axial specialization of $U_{\gl,\pp}$ to a nonsingular $M$-PAJ-$\gl$ factors through $U_{\gl,\pp}/\mathcal{R}(U_{\gl,\pp})$,  to get an axial specialization of $\overline{U_{\gl,\pp}}$
\end{proof}

\begin{remark}
    When $A$ is a PAJ-$\gl$, the radical of its normal Frobenius form corresponds with $\mathcal{R}_0,$ so the process stops at the first stage and the two concepts coincide.
\end{remark}

Now we go through the same steps as in Theorem~\ref{gen1}. By Theorem~\ref{TFa} we get finitely many generically $\mcT$-local generic MPAJ-$\gl$, $1\le \ell \le k',$ over integral domains, having the universal property that any nonsingular MPAJ-$\gl$ over $F$ is a marked homomorphic image of one of them.

 Write $\pi_\ell$ for the projection of $\overline{U_{\gl,\pp}}$ to ${\overline{U_{\gl,\pp}}}_\ell.$ 

\begin{thm}\label{gen1a}
 \begin{enumerate} \eroman
 \item Each $\overline{U_{\gl,N_\ell;M-\operatorname{prim}_\ell}},$ $1\le \ell\le k',$
is of Jordan type~$\gl$, and is generically $\mcT$-local.
 There is a canonical homomorphism $\Psi: \overline{U_{\gl,\pp}} \rightarrow \prod _{\ell =1}^{k'} \overline{U_{\gl,N_\ell;M-\operatorname{prim}_\ell}}$, whose kernel is $\mcT-\tor(\overline{U_{\gl,\pp}})$.
\item The assertions of Theorem~\ref{gen11} go over, mutatis mutandis, putting bars over all of the algebras.

 \end{enumerate}
\end{thm}
\begin{proof} The proofs of Theorems~\ref{gen1} and \ref{gen11} go over, mutatis mutandis.
\end{proof}

\subsubsection{The marked generic PAJ~$\gl$ with a normal Frobenius form}

We can incorporate
 a normal Frobenius form into our language, and use the Frobenius identity of \eqref{aya1} to build the generic PAJ~$\gl$ with a normal Frobenius form.

 \begin{lemma}\label{FF}
     Any Frobenius  form on an axial algebra $A$ generated by axes $\mcE$ is determined by $\{(a,y): a\in \mcE,\, y \text{ a word in } \mcE\}.$
 \end{lemma}
 \begin{proof}
     By bilinearity, we need only consider $(x,y)$ where $x,y$ are words in $\mcE$. We induct on the length $k$ of $x.$ Writing $x = za$ where $a\in \mcE$ and $z$ has length $k-1,$ $(x,y) = (za,y) = (z,ay) $, which satisfies the hypothesis by induction.
 \end{proof}
\begin{example}\label{Frb} Define $C^\sharp :=C^\natural [\chi_{i,j} : \ 1\le i\le N,\ 1\le j\le \theta(N,M)]$ and
$U_{\gl,N,M;\operatorname{preFrob}}=C^\sharp \{E\}$, also defining formally the symbols $(E_i,Y_j) = \chi_{i,j}$, with $\chi_{i,i}=1,$ $1\le i \le N.$ \footnote{The $\chi_{i,j}$ correspond somewhat to the $\gl_{[i,j]}$ of \cite[\S4]{HRS1}.}
To obtain generic objects, we need to quotient by the axial Frobenius identities of Remark~\ref{FI}, in view of Lemma~\ref{FF}.
\end{example}

 \begin{remark}
 This construction incidentally identifies $\chi_{i,j}$ with $\chi_{j,i}$ for $1\le i,j\le N$, since Frobenius forms are symmetric.

     The analogous construction would give us a generic axial algebra with an associative (not symmetric) bilinear form.
 \end{remark}

 \subsection{Matsuo algebras of Jordan type $\half$ versus Jordan algebras}

When $N\le 3$, we have seen that all PAJ-$\half$ are Jordan algebras. For $N\ge 4$, since there is a Matsuo algebra $A$ which is not Jordan even for Jordan type $\half$, cf.~\cite{DMR}, we could try to use it to reach a contradiction. But unfortunately the generic PAJ-$\half$ need not be torsion-free. So, we cannot answer Shpectorov's question at this stage.

\begin{remark}\label{pa1} In conclusion, for $N$ axially admissible, there may be several generic PAJ-$\half$ algebras over fields: several will be generic Jordan, several will be generic Matsuo, and the rest of them, \emph{if any exist}, will be neither Jordan nor the image of a Matsuo algebra, perhaps with some of them nonsingular.
In other words, a counterexample to Shpectorov's question would yield a generic counterexample.
\end{remark}

\section{Some further questions}

\begin{ques}
 \begin{enumerate}
 \item Can the nonsingular Baer radical defined here contain an axis $a$? (Not if $A$ has an anisotropic Frobenius form, since then the nonsingular Baer radical would coincide with the radical of the form.)

 \item How does the nonsingular Baer radical compare with  the radical defined in \cite[\S{3.2}]{GSh}, as the maximal ideal not containing an axis?

 \item What is the nonsingular Baer radical of $U_{\gl,\pp}$?
 \end{enumerate}
\end{ques}

\subsection{4-trace-admissibility}

 Long ago, Albert \cite{A3} introduced a certain kind of Frobenius form on power-associative algebras. We weaken it to algebras that need not be power-associative. But in this case, ``nilpotence'' of an element may be ambiguous; $x^2$ is unambiguous, as is $x^3 = x^2 x = x x^2$ when $A$ is commutative, but $x^4$ is ambiguous.

\begin{Def} Let $A$ be a commutative algebra.
\begin{enumerate} \eroman
\item An element $y\in A$ is \textbf{4-nilpotent} if $y^3 y =0$ and $y^2 y^2 = 0.$ A subset of $A$ is \textbf{4-nil} if each of its elements is 4-nilpotent.
 \item A normal Frobenius form on $A$ is \textbf{4-trace-admissible} if $(x,y)=0$ whenever $xy$ is 4-nilpotent.

 \item $A$ is \textbf{4-trace-admissible} if it has a 4-trace-admissible Frobenius form.

 \item A 4-trace-admissible algebra $A$ is \textbf{trace-admissible} if it also is power-associative.
\end{enumerate}
\end{Def}
\begin{remark}
 Albert \cite[p.~319]{A3} has an argument which proves that every power-associative algebra over a field, endowed with a Frobenius form for which no two nontrivial idempotents are orthogonal, is trace-admissible. However, his proof relied heavily on the fact that for any $y\in A$ the algebra $F[y]$ is associative and thus has a nontrivial idempotent unless $y$ is nilpotent.

 Schafer \cite{S} proved that every noncommutative Jordan algebra $A$ of characteristic~0 is trace-admissible.
 His proof passes to $A^+$, a Jordan algebra and thus trace-admissible, and trace-admissibility descends to $A.$
\end{remark}
\begin{ques}
\begin{enumerate}\eroman
 \item In a PAJ~$\half,$ is $(x,y) = 0$ whenever $xy$ is 4-nilpotent?
 \item More specifically, does this hold for Matsuo algebras?
 \end{enumerate}
\end{ques}

\bmhead{Acknowledgments}

The author thanks Yoav Segev for introducing him to this area and for many helpful discussions; Ivan Shestakov for important comments concerning almost Jordan and Lie triple algebras; Jari Desmet for detailed comments on earlier versions of the manuscript; and Sergey Shpectorov for his inspiring papers and influential work on axial algebras.

\section*{Statements and Declarations}

 \subsection{Funding} This work was supported by ISF grant 1994/20 and the Anshel Pfeffer Chair.
\subsection{Conflict of interest/Competing interests} None.

 \subsection{Ethics approval and consent to participate} Not applicable.
\subsubsection{Research involving Human Participants and/or Animals} None.

\subsubsection{Informed consent} Not applicable.

 \subsection{Consent for publication} Not applicable.
 \subsection{Data availability} Not applicable.
\subsection{Materials availability} Not applicable.
\subsection{Code availability} Not applicable.
\subsection{Author contributions} The author conceived and wrote the manuscript.

\subsection{Use of generative AI} AI was used solely to assist with checking statements and proofs, language editing, and bibliography. The author reviewed all suggested changes and takes full responsibility for the manuscript.

\end{document}